\documentstyle[11pt]{article}

\def\tabaddress#1{{\small\it\begin{tabular}[t]{c}#1 \\[1.2ex]\end{tabular}}}
\def\UPCMAT{Departamento de Matem\'atica Aplicada y Telem\'atica\\
   Universidad Polit\'ecnica de Catalu\~na\\
   Campus Nord, M\'odulo C-3\\
   C/ Gran Capit\'an s.n.\\
   E-08071 BARCELONA\\
   SPAIN}

\newtheorem{teor}{Theorem}
\newtheorem{prop}{Proposition}
\newtheorem{corol}{Corolary}
\newtheorem{lem}{Lemma}
\newtheorem{definition}{Definition}

\def\beq{\begin{equation}}
\def\eeq{\end{equation}}
\def\bea{\begin{eqnarray}}
\def\eea{\end{eqnarray}}
\def\beann{\begin{eqnarray*}}
\def\eeann{\end{eqnarray*}}
\def\beasn{\begin{sneqnarray}}
\def\eeasn{\end{sneqnarray}}
\def\ben{\begin{enumerate}}
\def\een{\end{enumerate}}
\def\bit{\begin{itemize}}
\def\eit{\end{itemize}}
\def\dst{\(\displaystyle}
\def\proof{( {\sl Proof} )\quad}
\def\derpar#1#2{\frac{\partial{#1}}{\partial{#2}}}

\def\map#1{\mathrel{\mathop{\to}\limits^{#1}}}
\def\mapping#1{\mathrel{\mathop{\longrightarrow}\limits^{#1}}}

\def\moment#1#2#3{{#1}_{#2}, \ldots, {#1}_{#3}}
\def\qed{\ifvmode\removelastskip\fi
{\unskip\nobreak\hfil\penalty50\hbox{}\nobreak\hfil
\hbox{\vrule height1.2ex width1.2ex}\parfillskip=0pt
\finalhyphendemerits=0 \par\smallskip}}

\def\vf{{\cal X}}
\def\df{{\mit\Omega}}
\def\Lag{{\cal L}}

\def\d{{\rm d}}

\def\Real{{\bf R}}
\def\Complex{{\bf C}}
\def\Ker{\mathop{\rm Ker}\nolimits}

\def\inn{\mathop{i}\nolimits}

\def\Tan{{\rm T}}

\def\Lie{\mathop{\rm L}\nolimits}
\def\Cinfty{{\rm C}^\infty}
\catcode`@=12
\def\ls{((E,M;\pi),\Lag )}

\def\lag{\pounds}
\def\del{{\cal E}^{\nabla}_{\Lag}}
\def\fel{{\rm E}^{\nabla}_{\Lag}}

\title{Geometry of Lagrangian First-order Classical Field Theories}
\author{\sc Arturo Echeverr\'ia-Enr\'iquez,
   \\
{\sc Miguel C. Mu\~noz-Lecanda\thanks{{\bf e}-{\it mail}: MATMCML@MAT.UPC.ES}},
   \\
{\sc Narciso Rom\'an-Roy\thanks{{\bf e}-{\it mail}: MATNRR@MAT.UPC.ES}},
   \\
   \tabaddress{\UPCMAT}}

\begin{document}
\maketitle
\thispagestyle{empty}

\begin{center}
{\bf DMAT 04-0194}
\end{center}

\begin{abstract}
We construct a lagrangian geometric formulation for
first-order field theories using the canonical structures of
first-order jet bundles, which are taken as the phase spaces of the systems
in consideration.
First of all, we construct all the geometric structures
associated with a first-order jet bundle and, using them,
we develop the lagrangian formalism,
defining the canonical forms associated with a lagrangian density
and the density of lagrangian energy, obtaining the
{\sl Euler-Lagrange equations} in two equivalent ways:
as the result of a variational problem and developing the
{\sl jet field formalism} (which is a formulation more similar
to the case of mechanical systems).
A statement and proof of Noether's theorem is also given,
using the latter formalism.

Finally, some classical examples are briefly studied.
\end{abstract}

\vfill \hfill
\vbox{\raggedleft AMS s.\,c.\,(1980): 53C80, 58A20. PACS: 0240, 0350 }\null

\clearpage

\tableofcontents

\section{Introduction}

Nowadays, it is well known that {\sl jet bundles} are the appropriate domain
for the description of classical field theory,
both in lagrangian and hamiltonian formalisms.
In recent years much work has been done in that domain
with the aim of establishing the suitable geometrical structures
for field theories in general
\cite{He-dgcv}, \cite{GS-73}, \cite{Gc-74}, \cite{AA-80},
\cite{GM-83}, \cite{LR-gcmft}, \cite{Gu-87}, \cite{Sa-87},
\cite{BSF-gcf}, \cite{Sa-89}, \cite{GIMMSY-mm}, \cite{CCI-91},
\cite{EM-91}, \cite{Sd-94b}
and for {\sl gauge} and {\sl Yang-Mills} field theories in particular
\cite{DM-fbt}, \cite{At-gymf}, \cite{Bl-gtvp}, \cite{DV-80},
\cite{EGH-ggt}, \cite{BK-ymtcs}, \cite{GS-dggt}, \cite{MM-ggf}.
Nevertheless, there is no definitive consensus about which are
the most relevant geometrical structures when dealing with
first-order lagrangian field theory and the role they play
\cite{Sa-89}, \cite{He-dgcv}, \cite{Gc-74},
\cite{GM-83}, \cite{Sa-87}, \cite{BSF-gcf}.
In addition, many authors have emphasized the importance of
the hamiltonian formulation of field theories
in the construction of the lagrangian formalism
\cite{Gu-87}, \cite{GIMMSY-mm}, \cite{CCI-91}, \cite{EM-91},
\cite{SO-92}, \cite{Sd-94a}
(in the same way that certain geometric formulations of mechanics
are made \cite{AM-78}).
However, even in these cases, an external element to the problem,
such as a bundle connection, is needed in order to develop
the hamiltonian field equations \cite{CCI-91}, \cite{EM-91}.

Our aim in this work is to construct a lagrangian geometric formulation for
first-order field theories using only the canonical structures of
the jet bundle $J^1E$ representing the system, in a way analogous to
{\sl Klein's formalism} for classical mechanics \cite{Kl-evm}.
A forthcoming paper will be devoted to a hamiltonian
formalism for these theories \cite{EMR-ghft}.

A brief description of the content of the work now follows:

\begin{itemize}
\item
The motivation of section 2 is the construction of the relevant
geometric structures of first-order jet bundles.
Following  \cite{Gc-74}, the {\sl vertical differential}
will be used as the fundamental tool in the build up of the
other elements: the {\sl canonical form}, the {\sl Cartan distribution},
and the {\sl vertical endomorphisms}.
We also study the {\sl canonical prolongations} of sections, diffeomorphisms
and vector fields and the invariance of the geometric elements by
these prolongations. Much of this section
is well known in the usual literature, but it is included here
in order to make this work more self-contained.
However, a new contribution is the construction of the vertical
endomorphisms (which play a most relevant role in the theory)
using the affine structure of the first-order jet bundle
and the canonical form. Another contribution is the proof of
the invariance of the geometric elements by the prolongation of
diffeomorphisms and vector fields above mentioned.
\item
In the first part of section 3 we develop the
lagrangian formalism. Using one of the vertical endomorphisms
we construct intrinsically the {\sl canonical forms}
associated with the {\sl lagrangian density}: the
{\sl Poincar\'e-Cartan forms}. A deep study is made about
the geometrical aspects of the variational problem posed by
a lagrangian density and we obtain the {\sl Euler-Lagrange equations},
both their local and global expressions, for sections of the
configuration bundle.

A particularly interesting point is the introduction
of the notion of the {\sl energy density} associated with a
lagrangian density and a connection. We give an intrinsic definition
of the energy density which, as expected, coincides with
the usual one when the configuration bundle is a direct product
and then a natural flat connection exists. Furthermore,
for a given lagrangian density we give the structure of the set
of all the possible energy densities associated with it.
We also establish the formula for the linear variation of the energy
density with respect to change on the connection.
We would like to point out that the need for a connection
for the definition of an energy density indicates a close
relationship, not developed in this work, between the energy
and the hamiltonian associated with the lagrangian and the connection
(see \cite{EMR-ghft}).

One of the main goals in this section is to establish
the {\sl Euler-Lagrange equations} for {\sl jet fields}
In order to sattisfy this aim, the problem of {\sl integrability}
and the {\sl second order condition} for jet fields are briefly treated.
Moreover, the solutions of such equations (if they exist)
are integrable second-order
jet fields whose integral manifolds are the critical sections
of the variational problem associated with the lagrangian.
In this sense, jet fields and their integral manifolds
play the same role as vector fields and their integral curves in mechanics.
So, in this language, field equations have a structural similarity with
the usual ones in mechanics.
We believe this formulation will, in due course,
provide us with a good opportunity
to tackle the study of {\sl singular lagrangian field theories}.
\item
A subject of great interest in field theory is the study of
{\sl symmetries}. Another section of the work is devoted to studying
the {\sl Noether symmetries} of a lagrangian problem
by means of the jet field formalism just developed.
We state and prove the {\sl Noether's theorem} in this formalism.
\item
A further section contains some typical examples, such as:
the {\sl electromagnetic field}, the {\sl bosonic string}
and the {\sl West-Zumino-Witten-Novikov model}.
For all of these we specify which is the corresponding
first-order jet bundle and its geometrical elements,
as well as the lagrangian density of the problem.
\item
The last section is devoted to discussion and outlook.
\item
The first appendix covers the main features on
the theory of connections on jet bundles and jet fields which
are used in some parts of the work.
Finally, since the notation on jet bundles is not unified,
a glosary of notation is also included at the end of the paper.
\end{itemize}

All the manifolds are real, second countable and $\Cinfty$.
The maps are assumed to be $\Cinfty$.
Sum over repeated indices is understood.

\section{Elements of differential geometry for
first-order lagrangian field theories}

\subsection{Geometrical structures of first-order jet bundles}

In this first paragraph we introduce the basic definitions and
properties concerning to first-order jet bundles, as well as the
canonical structures defined therein.
For more details regarding some of them see \cite{Sa-89}.

\subsubsection{First-order jet bundles}
\protect\label{fojb}

Let $M$ be an orientable manifold and $\pi\colon E\longrightarrow M$
a differentiable fiber bundle with typical fibre $F$. We denote by
$\Gamma (M,E)$ or $\Gamma (\pi)$ the set of global sections
of $\pi$. In the same way, if $U\subset M$ is an open set,
let $\Gamma_U(\pi )$ be the set of local sections of $\pi$
defined on $U$. Let $\dim M=n+1$, $\dim F=N$.

{\bf Remarks}:
\begin{itemize}
\item
The dimension of $M$ is taken to be $n+1$ because,
in many cases, this manifold is space-time. The bundle
$\pi\colon E\longrightarrow M$ is called the
{\sl covariant configuration bundle}.
The {\sl physical fields} are the sections of this bundle.
\item
The orientability of $M$ is not used in this first part
of the work. It is relevant when dealing with variational problems.
\end{itemize}

We denote by $J^1E$ the {\sl bundle of $1$-jets} of sections of $\pi$,
or {\sl $1$-jet bundle}, which is
endowed with the natural projection $\pi^1 \colon J^1E\longrightarrow E$.
For every $y\in E$, the fibers of $J^1E$ are denoted $J^1_yE$
and their elements by $\bar y$.
If $\phi \colon U\rightarrow E$ is a representative of
$\bar y \in J^1_yE$, we write $\phi \in \bar y$
or $\bar y=\Tan_{\pi (y)}\phi$.

In addition,
the map $\bar\pi^1 = \pi \circ \pi^1\colon J^1E \longrightarrow M$
defines another structure of differentiable bundle.
We denote by ${\rm V}(\pi)$ the vertical bundle associated with $\pi$,
that is ${\rm V}(\pi)=\Ker\Tan\pi$, and by
${\rm V}(\pi^1 )$ the vertical bundle associated with $\pi^1$,
that is ${\rm V}(\pi^1)=\Ker\Tan\pi^1$.
$\vf^{{\rm V}(\pi)}(E)$ and $\vf^{{\rm V}(\pi^1 )}(J^1E)$
will denote the corresponding sections or vertical
vector fields. Then, it can be proved (see \cite{Ou-73}) that
$\pi^1 \colon J^1E\longrightarrow E$
is an affine bundle modelled on the vector bundle
${\bf E}=\pi^*\Tan^*M \otimes_E{\rm V}(\pi)$
\footnote{
This notation denotes the tensor product of two vector bundles
over $E$.
}.
Therefore, the rank  of $\pi^1\colon J^1E\longrightarrow E$ is $(n+1)N$.

Sections of $\pi$ can be lifted to $J^1E$ in the following way:
let $\phi \colon U\subset M \to E$ be a local section of $\pi$.
For every $x\in U$, the section $\phi$ defines  an element of
$J^1E$: the equivalence class of $\phi$ in $x$, which is denoted
$(j^1\phi )(x)$.
Therefore we can define a local section of $\bar\pi^1$
$$
\begin{array}{ccccc}
j^1\phi &\colon &U & \longrightarrow & J^1E \\
& & x & \longmapsto     & (j^1 \phi )(x)
\end{array}
$$
and so we have defined a map
$$
\begin{array}{ccc}
\Gamma_U(M,E) & \longrightarrow & \Gamma_U(M,J^1E)  \\
       \phi         & \longmapsto     & j^1(\phi )\equiv j^1\phi
\end{array}
$$
$j^1\phi$ is called the {\rm canonical lifting} or the
{\sl canonical prolongation} of $\phi$ to $J^1E$.
A section of $\bar\pi^1$ which is the canonical extension of a section
of $\pi$ is called a {\sl holonomic section}.

Let $x^\mu$, $\mu = 1,\ldots,n,0$, be a local system in $M$
and $y^A$, $A= 1,\ldots,N$ a local system in the fibers;
that is, $\{ x^\mu ,y^A\}$ is a coordinate system adapted to the bundle.
In these coordinates, a local section
$\phi\colon U \rightarrow E$ is writen as
$\phi (x)=(x^\mu,\phi^A(x))$, that is, $\phi (x)$
is given by functions $y^A=\phi^A(x)$.
These local systems ${x^\mu}$, ${y^A}$ allows us to construct a
local system $(x^\mu,y^A,v^A_\mu)$ in $J^1E$, where
$v^A_\mu$ are defined as follows: if $\bar y\in J^1E$,
with $\pi^1 (\bar y)=y$ and $\pi(y)=x$,
let $\phi\colon U\rightarrow E$, $y^A=\phi ^A$, be a representative
of $\bar y$, then
$$
v^A_\mu (\bar y)=\left(\derpar{\phi^A}{x^\mu}\right)_x
$$
These systems are called {\sl natural local systems} in $J^1E$.
In one of them we have
$$
j^1\phi (x)=\left(x^\mu (x),\phi ^A (x),\derpar{\phi^A}{x^\mu}(x)\right)
$$

{\bf Remarks}:
\begin{itemize}
\item
If $\bar y\in J^1_yE$, with $x=\pi(y)$, and $\phi\colon U\to E$
is a representative of $\bar{y}$ we have the split
$$
\Tan_yE={\rm Im}\ \Tan_x\phi\oplus{\rm V}_y(\pi)
$$
hence the sections of $\pi^1$ are identified with connections
in the bundle $\pi\colon E\longrightarrow M$,
since they induce a horizontal subbundle of $\Tan E$.

Observe that it is reasonable to write ${\rm Im}\ \bar y$ for an
element $\bar y\in J^1E$.
\item
Giving a global section of an affine bundle,
this one can be identified with its associated vector bundle.
In our case we have:
\begin{itemize}
\item
If $\pi\colon E\longrightarrow M$ is a trivial bundle, that is
$E=M\times F$, a section of $\pi^1$ can be chosen in the following way:
denoting by $\pi_1\colon M\times F\to M$ and $\pi_2\colon M\times F\to F$
the canonical projections,
for a given $y_o\in M\times F$, $y_o=(x_o,v_o)=(\pi_1(y_o),\pi_2(y_o))$,
we define the section
$\phi_{y_o}(x)=(x,\pi_2(y_o))$, for every $x\in M$.
With these sections of $\pi$ we construct another one
$z$ of $\pi^1$ as follows:
$$
z(y):=(j^1\phi_y)(\pi_1(y))\quad ; \quad y\in E
$$
which is taken as the zero section of $\pi^1$.
So, in this case, $J^1E$ is a vector bundle over $E$.
\item
If $\pi\colon E\longrightarrow M$ is a vector bundle with
typical fiber $F$, let $\phi :M\to E$ be
the zero section of $\pi$ and $j^1\phi\colon M\to J^1E$
its canonical lifting. We construct the zero section of
$\bar\pi^1$ in the following way:
$$
z(y):=(j^1\phi)(\pi(y))\quad ; \quad y\in E
$$
thereby, in this case $\pi^1\colon J^1E\longrightarrow E$
is a vector bundle.
\end{itemize}
\end{itemize}

\subsubsection{Vertical differential. Canonical form}

(Following \cite{Gc-74}).

We are now going to define the canonical geometric structures with which
$J^1E$ is endowed with. The first one is the
{\sl canonical form}. Firstly we need to introduce
the concept of {\sl vertical differentiation}:

\begin{definition}
Let $\phi\colon M\to E$ be a section of $\pi$, $x\in M$ and $y=\phi(x)$.
The {\rm vertical differential} of the section $\phi$ at the point $y\in E$
is the map
$$
\begin{array}{ccccc}
\d^v_y\phi&\colon&\Tan_y E & \longrightarrow & {\rm V}_y(\pi) \\
  & & u & \longmapsto & u-\Tan_y(\phi\circ\pi)u
\end{array}
$$
\end{definition}

Notice that $\d^v_y \phi$ is well defined since
$\Tan_y \pi \circ \d^v_y \phi = 0$,
so the image is in ${\rm V}_y(\pi)$ and it depends only on $(j^1 \phi)(x)$.

If $(x^\mu ,y^A)$ is a natural local system of
$E$ and $\phi=(x^\mu ,\phi^A(x^\mu ))$, then
$$
\d^v_y \phi \left(\derpar{}{x^\mu}\right)_y =
-\left(\derpar{\phi^A}{x^\mu}\right)_y\left(\derpar{}{y^A}\right)_y
\quad , \quad
\d^v_y \phi \left(\derpar{}{y^A}\right)_y =
\left(\derpar{}{y^A}\right)_y
$$

Observe that the vertical differential splits $\Tan_y E$
into a vertical component and another one which is tangent
to the imagen of $\phi$ at the point $y$.
In addition, $\d^v_y \phi (u)$ is the projection of $u$
on the vertical part of this splitting.

As $\d^v_y \phi$ depends only on $(j^1\phi )(\pi (y))$,
the vertical differential can be lifted to $J^1E$
in the following way:

\begin{definition}
Consider $\bar y\in J^1E$ with
$\bar y\stackrel{\pi^1}{\mapsto}y\stackrel{\pi}{\mapsto}x$
and $\bar u\in\Tan_{\bar y}J^1E$.
The {\rm structure canonical form} of $J^1E$ is a $1$-form $\theta$
in $J^1E$ with values on ${\rm V}(\pi )$
which is defined by
$$
\theta (\bar y;\bar u):=(\d^v_y \phi)(\Tan_{\bar y}\pi^1 (\bar u))
$$
where the section $\phi$ is a representative of $\bar y$.
\end{definition}

This expression is well defined and does not depend on
the representative $\phi$ of $\bar y$.

If $(x^\mu ,y^A,v^A_\mu )$ is a natural local system of
$J^1E$ and $x^\mu (\bar y)=\alpha^\mu$,
$y^A(\bar y)=\beta^A$, $v^A_\mu (\bar y)=\gamma^A_\mu$, then:
$$
\theta \left(\bar y;\left(\derpar{}{x^\mu}\right)_{\bar y}\right) =
-\gamma^A_\mu\left(\derpar{}{y^A}\right)_y
\quad , \quad
\theta\left(\bar y;\left(\derpar{}{y^A}\right)_{\bar y}\right) =
\left(\derpar{}{y^A}\right)_y
\quad , \quad
\theta\left(\bar y;\left(\derpar{}{v^A_\mu}\right)_{\bar y}\right) = 0
$$
where $\phi$ is a representative of $\bar y$.
As one may observe, the images can be placed in the point
$\bar y$ lifting them to this point from $y$,
by means of the pull-back of bundles.

 From the above calculations we conclude that $\theta$ is
diferentiable and its expression in a natural local system is
\beq
\theta =  ({\rm d}y^A - v^A_{\mu}{\rm d}x^{\mu})
\otimes \derpar{}{y^A}
\label{canform}
\eeq
According to this, $\theta$ is  a differential $1$-form
defined in $J^1E$ with values on $\pi^{1^*}{\rm V}(\pi)$,
hence it is an element of
$\df^1(J^1E,\pi^{1^*}{\rm V}(\pi))=
\df^1(J^1E)\otimes_{J^1{E}}\Gamma (J^1E,\pi^{1^*}{\rm V}(\pi))$,
where $\df^1(J^1E)$ denotes the set of differential $1$-forms in $J^1E$.

Holonomic sections can be characterized using the
structure canonical form as follows:

\begin{prop}
Let $\psi\colon M \to J^1E$ be a section of $\bar\pi^1$.
The necessary and sufficient condition for $\psi$
to be a holonomic section is that
$\psi^*\theta=0$
\end{prop}
( {\sl Proof} )\quad
($\Longleftarrow$)\quad
Consider $\psi =j^1\phi$. If $x\in M$ and $v\in\Tan_xM$, we have
\begin{eqnarray*}
(\psi^*\theta)(x;v) &=& ((j^1\phi)^*\theta)(x;v)
=\theta((j^1\phi)(x);(\Tan_xj^1\phi)v)
\\ &=&
(\d^v_{\phi(x)}\phi)((\Tan_{(j^1\phi)(x)}\pi^1 )(\Tan_x j^1\phi(v))) =
(\d^v_{\phi(x)}\phi)((\Tan_x \phi)(v))
\\ &=&
({\rm T}_x \phi)(v)
-{\rm T}_{\phi(x)}(\phi\circ \pi)(({\rm T}_x \phi)(v))=0
\end{eqnarray*}

\quad\quad ($\Longrightarrow$)\quad
Now, let $\psi\colon M\to J^1E$ be a section of $\bar\pi^1$
such that $\psi^*\theta =0$. If $x\in M$ and $(x^\mu,y^A,v^A_\mu )$ is a
natural system of coordinates in a neighbourhood $U$ of $\psi(x)$;
then we have $\psi\vert_U=(x^\mu,g^A(x^\mu),f^A_\rho(x^\mu))$
and therefore
$$
0=\psi^*\theta =\left(\derpar{g^A}{x^\mu}dx^\mu
-f^A_\mu dx^\mu\right)\otimes \derpar{}{y^A}
$$
then \dst f^A_\mu(x^\rho)=\derpar{g^A}{x^\mu}(x^\rho)\) and $\psi$ is the
canonical prolongation of a section of $\pi$.
\qed

We will use this canonical form  in order to construct some
other geometric elements on $J^1E$.

\subsubsection{Contact module (Cartan distribution)}

The second canonical structure of $J^1E$ is defined in the following way:

\begin{definition}
Let $\theta$ be the structure canonical form of $J^1E$.
As we have seen, it is an element of
$\df^1(J^1E)\otimes_{J^1E}\Gamma (J^1E,\pi^{1^*}{\rm V}(\pi))$;
then it can be considered as a $\Cinfty (J^1E)$-linear map
$$
\theta\colon\Gamma (J^1E,\pi^{1^*}{\rm V}(\pi))^*\longrightarrow\df^1(J^1E)
$$
The image of this map is called the {\rm contact module}
or {\rm Cartan distribution} of $J^1E$.
It is denoted by ${\cal M}_c$.
\end{definition}

\begin{prop}
${\cal M}_c$ is a locally finite generated module.
If $(x^\mu ,y^A,v^A_\mu )$ is a local natural system of coordinates
in an open set of $J^1E$, then ${\cal M}_c$ is generated by the forms
$\theta^A=\d y^A-v^A_\mu\d x^\mu$ in this open set.
\label{finmod}
\end{prop}
( {\sl Proof} )\quad
The local expression of $\theta$ in this local system is
\dst\theta =\left(\d y^A-v^A_\mu\d x^\mu\right)\otimes\derpar{}{y^A}\) .
Let $\{\xi^A\}$ be the local basis of $\Gamma (J^1E,\pi^{1^*}{\rm V}(\pi))^*$
which is dual of \dst\left\{\derpar{}{y^A}\right\}\) ,
therefore the image of $\theta$ has
as a local basis the above mentioned forms.
\qed

{\bf Remark}:
\bit
\item
Observe that there is no identification of ${\rm V}(\pi)^*$
as a subbundle of $\Tan^*E$ unless we have a connection on
$\pi\colon E\to M$.
So the elements $\xi^A$ are defined by
\dst\xi^A\left(\derpar{}{y^B}\right) =\delta^A_B\) .
\eit

Now, using the Cartan distribution, we have another characterization
of holonomic sections:

\begin{prop}
Let $\psi\colon M \to J^1E$ be a section of $\bar\pi^1$.
The necessary and sufficient condition for
$\psi$ to be a holonomic section is that
$\psi^*\alpha =0$, for every $\alpha\in {\cal M}_c$.
\end{prop}
( {\sl Proof} )\quad
($\Longrightarrow$)\quad
Let $\psi=j^1\phi$ be an holonomic section and $\alpha\in {\cal M}_c$.
By definition $\alpha =\theta (\varphi )$, where
$\varphi\colon J^1E\to\pi^{1^*}{\rm V}(\pi)$ is a section.
If $x\in M$ and $v\in\Tan_xM$, we have
\beann
(\psi^*\alpha )(x;v)&=&
\alpha (\psi(x);(\Tan_x\psi)v)=\theta (\varphi )(\psi(x);(\Tan_x\psi)v)
\\&=&
\theta (\psi(x);(\Tan_x\psi)v,\varphi )=
\theta (\psi(x);(\Tan_x\psi)v)(\varphi )=
((\psi^*\theta )(x;v)(\varphi )=0
\eeann
since $\psi^*\theta =0$, because $\psi$ is holonomic.

\quad\quad ($\Longleftarrow$)\quad
It suffices to see it in an open set of coordinates.
Let $(x^\mu ,y^A,v^A_\mu )$ be a local natural system of coordinates
in $(\bar\pi^1)^{-1}(U)$, with $U$ an open set in $M$.
Let $\psi\colon M\to J^1E$ be a section such that $\psi^*\alpha=0$,
for every $\alpha\in {\cal M}_c$. In particular, in $U$, we have
$\psi\vert_U=(x^\mu ,g^A(x^\mu ),f^B_\rho(x^\mu ))$,
and $\psi$ vanishes on $\theta^A$,
hence \dst f^B_\rho =\derpar{g^B}{x^\rho}$ and $\psi$ is holonomic.
\qed

Conversely, the forms belonging to ${\cal M}_c$ can be characterized
by means of holonomic sections, which is the classical
definition of ${\cal M}_c$:

\begin{prop}
Consider $\alpha\in\df^1(J^1E)$.
The necessary and sufficient condition for
$\alpha\in {\cal M}_c$ is that $(j^1\phi)^*\alpha =0$,
for every section $\phi\colon M\to E$.
\end{prop}
( {\sl Proof} )\quad
We define the set
$$
{\cal M}_0 := \{\alpha\in\df^1(J^1E)\ |\ (j^1\phi )^*\alpha =0,\
\forall \phi\colon M\to E, {\rm section\ of}\ \pi\}
$$
It is clear that ${\cal M}_c\subset {\cal M}_0$. Now we are able to see that
${\cal M}_0$ coincides with ${\cal M}_c$ locally.

Consider $\bar y\in J^1E$ with $\pi^1 (\bar y)=y$ and
$\pi(y)=x$, and let $\phi\colon M\to E$ be a section of
$\pi$ with $\phi(x)=y$ and $(j^1\phi)(x)=\bar y$.
If $(x^\mu ,y^A,v^A_\mu ))$ is a local system of coordinates and
$\phi=(x^\mu ,\phi^A(x^\mu )$ in this local system, we have
\beann
(\Tan_xj^1\phi)\derpar{}{x^\mu}\Big\vert_x &=&
\derpar{}{x^\mu}\Big\vert_{\bar y}+
\derpar{\phi^A}{x^\mu}\Big\vert_x\derpar{}{y^A}\Big\vert_{\bar y}+
\frac{\partial^2\phi^A}{\partial x^\mu\partial x^\rho}\Big\vert_x
\derpar{}{v^A_\rho}\Big\vert_{\bar y}
\\ &=&
\derpar{}{x^\mu}\Big\vert_{\bar y}+
v^A_\mu (\bar y)\derpar{}{y^A}\Big\vert_{\bar y}+
\frac{\partial^2\phi^A}{\partial x^\mu\partial x^\rho}\Big\vert_x
\derpar{}{v^A_\rho}\Big\vert_{\bar y}
\eeann
Therefore, varying the section $\phi$, we obtain all the vectors in
$\Tan_{\bar y}J^1E$ which are tangent to the image of the canonical extension
of sections of $\pi$. They are all the vectors of the subspace
generated by
$$
\derpar{}{x^\mu}\Big\vert_{\bar y}+
v^A_\mu (\bar y)\derpar{}{y^A}\Big\vert_{\bar y}
\quad ; \quad
\derpar{}{v^A_\rho}\Big\vert_{\bar y}
$$
A basis of the subspace of $\Tan^*_{\bar y}J^1E$ which
is incident to the last one is
$$
\theta^A_{\bar y}=(\d y^A-v^A_\mu (\bar y)\d x^\mu )\vert_{\bar y}
$$
thereby ${\cal M}_0$ is locally generated by
$\theta^A=\d y^A-v^A_\mu\d x^\mu$ and, hence, it coincides with ${\cal M}_c$.
\qed

\subsubsection{Vertical endomorphisms}

Finally, we have two geometrical elements which generalize
the notion of the canonical endomorphism of a tangent bundle
\cite{Kl-evm}.

As we have seen, $J^1(E)$ is an affine bundle over $E$ with associated
vector bundle $\pi^*\Tan^*M\otimes_E{\rm V}(\pi)$.
Therefore, if $\bar y \in J^1_yE$, the tangent space to the fiber
$J^1_yE$ at $\bar y$ is canonically isomorphic to
$\Tan_{\pi(y)}^*M\otimes{\rm V}_y(\pi)$.
Moreover, $\Tan_{\bar y}J^1_yE$ is just ${\rm V}_{\bar y}J^1E$,
that is, the vertical tangent space of $J^1E$
with respect to the projection $\pi^1$.
Then:

\begin{definition}
For every $\bar y\in J^1E$, consider the canonical isomorphism
$$
{\cal S}_{\bar y}\colon
\Tan_{\bar\pi^1 (\bar y)}^*M\otimes{\rm V}_{\pi^1 (\bar y)}(\pi)
\longrightarrow {\rm V}_{\bar y}(\pi^1 )
$$
which consists in associating with an element
$\alpha\otimes v\in
\Tan_{\bar\pi^1 (\bar y)}^*M\otimes{\rm V}_{\pi^1 (\bar y)}(\pi)$
the directional derivative in $\bar y$ with respect to $\alpha\otimes v$.
Taking into account that  $\alpha\otimes v$ acts in $J^1_yE$
by translation
$$
{\cal S}_{\bar y}(\alpha\otimes v ):= D_{\alpha\otimes v}(\bar y)
\colon f\mapsto \lim_{t\to 0}\frac{f(\bar y+t(\alpha\otimes v))-
f(\bar y)}{t}
$$
for $f\in\Cinfty (J^1_yE)$.
Then we have the following isomorphism of $\Cinfty (J^1E)$-modules
$$
{\cal S}\colon
\Gamma (J^1E,\bar\pi^{1^*}\Tan^*M)\otimes\Gamma (J^1E,\pi^{1^*}{\rm V}(\pi))
\longrightarrow \Gamma (J^1E,{\rm V}(\pi^1 ))
$$
which is called the {\rm vertical endomorphism} ${\cal S}$.
\end{definition}

If $(x^\mu ,y^A,v^A_\mu )$ is a natural system of coordinates,
the local expression of ${\cal S}$ is given, in this local system, by
$$
{\cal S}\left(\d x^\mu\otimes\derpar{}{y^A}\right)=\derpar{}{v^A_\mu}
$$
and, therefore,
$$
{\cal S}=\xi^A\otimes\derpar{}{v^A_\mu}\otimes\derpar{}{x^\mu}
$$
where $\{\xi^A\}$ is the local basis of $\Gamma (J^1E,\pi^{1^*}{\rm V}(\pi))^*$
which is dual of \dst\left\{\derpar{}{y^A}\right\}\)
as we have pointed out in the remark following proposition
\ref{finmod}. So we have:

\begin{prop}
${\cal S}\in\Gamma (J^1E,(\pi^{1^*}{\rm V}(\pi))^*)\otimes
\Gamma (J^1E,{\rm V}(\pi^1 ))\otimes\Gamma(J^1E,\bar\pi^{1^*}\Tan M)$
(where all the tensor products are on $\Cinfty (J^1E)$).

In addition, ${\cal S}$ is a canonical section of
$(\pi^{1^*}{\rm V}(\pi))^*\otimes{\rm V}(\pi^1 )\otimes\bar\pi^{1^*}\Tan M$
over $J^1E$ (canonical in the sense that it is associated with the
initial structure $\pi\colon E\to M$).
\end{prop}

Moreover, we have the canonical form
$\theta\in\df^1(J^1E)\otimes_{J^1E}\Gamma (J^1E,\pi^{1^*}{\rm V}(\pi))$,
then we can define:

\begin{definition}
The {\rm vertical endomorphism} ${\cal V}$ arises from
the natural contraction between
the factors $\Gamma (J^1E,(\pi^{1^*}{\rm V}(\pi))^*)$ of ${\cal S}$
and $\Gamma (J^1E,\pi^{1^*}{\rm V}(\pi))$ of $\theta$; that is
$$
{\cal V}=
\inn ({\cal S})\theta\in\df^1(J^1E)\otimes\Gamma (J^1E,{\rm V}(\pi^1 ))
\otimes\Gamma (J^1E,\bar\pi^{1^*}\Tan M)
$$
\end{definition}

Taking into account the local expression of $\theta$ (eq. (\ref{canform})),
it is clear that
the local expression of ${\cal V}$ in a natural system of coordinates is
$$
{\cal V}=\left(\d y^A-v^A_\mu\d x^\mu\right)\otimes
\derpar{}{v^A_\nu}\otimes\derpar{}{x^\nu}
$$

$J^1E$ is also endowed with other canonical structures, namely:
the {\sl total differentiation} and the
{\sl module of total derivatives} (see \cite{Sa-89}).
Nevertheless we are not going to use them in this work.

\subsection{Canonical prolongations}
\protect\label{prol}

Consider $J^1E\mapping{\pi^1}E\mapping{\pi}M$ as above.
Geometric objects in $E$ (namely, differential forms and vector fields)
can be canonically lifted to $J^1E$.
Next, we study the way this is done, as well as the main properties
of these liftings, because they will be used later
in order to characterize critical sections starting from
variational principles and in the study of symmetries.

\subsubsection{Prolongation of sections and diffeomorphisms}

Let $\phi\colon U\subset M\to E$ be a local section.
In the above paragraphs we have seen that this section can be prolonged
to a section $j^1\phi\colon U\to J^1E$ with the following properties
$$
1. \quad \pi^1\circ j^1\phi=\phi
\qquad\qquad
2. \quad (j^1\phi)^*\theta =0
$$
In order to evaluate $j^1\phi$ in a point $x\in U$, it suffices to take
the equivalence class of $\phi$ in $x$, as we have described in the
construction of $J^1E$.

Now, let $\Phi\colon E\to E$ be a diffeomorphism of
$\pi$-fiber bundles and $\Phi_M\colon M\to M$ the diffeomorphism
induced on the basis. Our aim is to lift $\Phi$ to a diffeomorphism
of $J^1E$ in a natural way; that is, in such a way that the following
diagram commutes:
$$
\begin{array}{rcccl}
& J^1E & \mapping{j^1\Phi} & J^1E & \\
\pi^1 & \Big\downarrow & & \Big\downarrow & \pi^1 \\
& E & \mapping{\Phi} & E & \\
\pi & \Big\downarrow & & \Big\downarrow & \pi \\
& M & \mapping{\Phi_M} & M & \\
\end{array}
$$

\begin{definition}
Consider $\bar y\in J^1E$, with
$\bar y\stackrel{\pi^1}{\mapsto}y\stackrel{\pi}{\mapsto}x$,
and let $\phi\colon U\to E$ be a local representative of $\bar y$; that is
$\phi\in\bar y$. Then $\Phi\circ\phi\circ\Phi_M^{-1}$
is a local section of $\pi$ defined on $\Phi_M(U)$.
The {\rm canonical prolongation} or the {\rm canonical lifting}
of the diffeomorphism $\Phi$ is
the map $j^1\Phi\colon J^1E\longrightarrow J^1E$
defined by
$$
(j^1\Phi)(\bar y):=j^1(\Phi\circ \phi\circ\Phi_M^{-1})(\Phi_M(x))
$$
It is clear that this definition does not depend on the
representative $\phi\in\bar y$; that is, $(j^1\Phi)(\bar y)$ is well
defined and, moreover, it depends differentiably on $\bar y$.
\end{definition}

We can check that this definition has the predicted properties. In fact:

\begin{prop}
Let $\Phi\colon E\to E$ be a diffeomorphism of fiber bundles. Then
\ben
\item
$\pi^1\circ j^1\Phi =\Phi\circ\pi^1$,
$\bar\pi^1\circ j^1\Phi =\Phi_M\circ\bar\pi^1$.
\item
If $\Phi'\colon E\to E$ is another fiber bundle diffeomorphism, then
$$
j^1(\Phi'\circ\Phi_1)=j^1\Phi'\circ j^1\Phi
$$
\item
$j^1\Phi$ is a diffeomorphism of $\pi^1$-bundles and $\bar\pi^1$-bundles, and
$(j^1\Phi )^{-1}=j^1\Phi^{-1}$.
\item
If $\phi\colon U\to E$ is a local section of $\pi$, then
$$
j^1(\Phi\circ \phi\circ\Phi_M^{-1})=j^1\Phi\circ j^1\phi\circ\Phi_M^{-1}
$$
\item
If $\psi\colon U\to J^1E$ is a holonomic section,
with $\psi =j^1\phi$ for some $\phi\colon U\to E$,
then $j^1\Phi\circ j^1 \phi\circ\Phi_M^{-1}\colon\Phi_M(U)\to J^1E$
is a holonomic section.
\een
\label{ppc}
\end{prop}
( {\sl Proof} )
See \cite{Sa-89}.
\qed

{\bf Remark}:
\begin{itemize}
\item
If $\Phi\colon E\to E$ induces the identity in $M$
(that is a {\sl strong diffeomorphism}), then $j^1\Phi$
is defined by
$$
(j^1\Phi )(\bar y)=j^1(\Phi\circ \phi)(x)
\quad ; \quad (\phi\in\bar y \quad , \quad  \bar\pi^1 (\bar y)=x)
$$
since $\Phi\circ \phi$ is a section of $\pi$ defined in the same
open set as the section $\phi$.
In the same way, in this case, $j^1\Phi\circ j^1 \phi=j^1(\Phi\circ \phi)$,
and hence canonical prolongations of strong diffeomorphisms
transform holonomic sections
defined in an open set into holonomic sections defined in the same open set.
\end{itemize}

\subsubsection{Relation with the affine structure}
\protect\label{ras}

Now we want to describe the relations between the canonical prolongations of
diffeomorphisms and the canonical structures of $J^1E$ introduced in
the above section.

First of all we study the relation with the affine structure of $J^1E$.

\begin{prop}
The canonical prolongations
of diffeomorphisms $\Phi\colon E\to E$ over $\pi$ preserve
the affine structure of $\pi^1\colon J^1E\to E$.
\label{canprol}
\end{prop}
( {\sl Proof} )\quad
Consider $\bar y\in J^1E$ and $\phi\in\bar y$, we have
$$
(j^1\Phi )(\bar y)=
(j^1(\Phi\circ \phi\circ\Phi_M^{-1})\circ\Phi_M)(\bar\pi^1(\bar y))
$$
Consider now
$\alpha\otimes v\in
\Tan^*_{\bar\pi^1 (y)}M\otimes{\rm V}_{\pi^1 (\bar y)}(\pi)$
and the point $\bar y+\alpha\otimes v$.
Let $\phi'\colon M\to E$ be a local section with
$\phi'(\bar\pi^1 (\bar y))=\pi^1 (\bar y)$ and
$\Tan_{\bar\pi^1 (\bar y)}\phi'=\Tan_{\bar\pi^1 (\bar y)}\phi+\alpha\otimes v$.
We have
\beann
(j^1\Phi )(\bar y+\alpha\otimes v)&=&
(j^1(\Phi\circ \phi'\circ\Phi_M^{-1})\circ\Phi_M)
(\bar\pi^1 (\bar y+\alpha\otimes v))
\\ &=&
(j^1(\Phi\circ \phi'\circ\Phi_M^{-1})\circ\Phi_M)(\bar\pi^1 (\bar y))
\\ &=&
\Tan_{\phi'(\bar\pi^1 (\bar y))}\Phi\circ\Tan_{\bar\pi^1 (\bar y)}\phi'
\circ\Tan_{\Phi_M(\bar\pi^1 (\bar y))}\Phi_M^{-1}
\\ &=&
\Tan_{\phi(\bar\pi^1 (\bar y))}\Phi\circ (\Tan_{\bar\pi^1
(\bar y)}\phi+\alpha\otimes v)\circ\Tan_{\Phi_M(\bar\pi^1 (\bar y))}\Phi_M^{-1}
\\ &=&
\Tan_{\phi(\bar\pi^1 (\bar y))}\Phi\circ\Tan_{\bar\pi^1 (\bar y)}\phi
\circ\Tan_{\Phi_M(\bar\pi^1 (\bar y))}\Phi_M^{-1} +
\Tan_{\phi(\bar\pi^1 (\bar y))}\Phi\circ\alpha\otimes v\circ
\Tan_{\Phi_M(\bar\pi^1 (\bar y))}\Phi_M^{-1}
\\ &=&
j^1(\Phi\circ \phi\circ\Phi_M^{-1})(\Phi_M(\bar\pi^1 (\bar y))+
(\Tan_{\Phi_M(\bar\pi^1 (\bar y))}^*\Phi_M^{-1}\otimes\Tan_{\pi^1 (\bar y)})
\Phi (\alpha\otimes v)
\\ &=&
(j^1\Phi )(\bar y)+\varphi (\alpha\otimes v)
\eeann
where
$\varphi\colon\Tan_{\bar\pi^1 (\bar y)}^*M\otimes{\rm V}_{\pi^1 (\bar y)}(\pi )
\to \Tan_{\Phi_M(\bar\pi^1 (\bar y))}^*M\otimes
{\rm V}_{\Phi (\bar\pi^1 (\bar y))}(\pi )$
is the map
$\varphi:=
\Tan_{\Phi_M(\bar\pi^1 (\bar y))}^*\Phi_M^{-1}\otimes\Tan_{\pi^1 (\bar y)}\Phi$
and hence it is linear, therefore $j^1\Phi$ is an affine mapping.
\qed

{\bf Remarks}:
\begin{itemize}
\item
See \cite{Sa-89} for a local version of this proof.
\item
Consider $\bar y\in J^1E$ with $\pi^1 (\bar y)=y$.
Since $J^1_yE$ is an affine bundle over $E$ with typical fiber
$\bar\pi^{1^*}\Tan^*M\otimes\pi^{1^*}{\rm V}(\pi)$,
the tangent space to one fiber at
each point is canonically isomorphic to the vector space where
the fiber is modelled on; that is
$$
\Tan_{\bar y}^*J^1_yE\simeq\Tan_{\bar\pi^1 (\bar y)}^*M\otimes{\rm V}_y(\pi)
$$
If $\Phi\colon E\to E$
is a $\pi$-diffeomorphism, we have seen that $j^1\Phi$
preserves the affine structure, hence the tangent of $j^1\Phi$
is just its linear part; that is, $\Tan_{\bar y}j^1\Phi$
acts on $\Tan_{\bar y}J^1E$ in the following way:
$$
\begin{array}{ccc}
\Tan_{\bar y}J^1E& $\rightarrowfill$ &\Tan_{j^1\Phi (\bar y)}J^1_{\Phi (y)}E
\\
&\Tan_{\bar y}j^1\Phi&
\\
\simeq\Big\updownarrow & &\Big\updownarrow\simeq
\\
&\Tan_{\Phi_M(\bar\pi^1 (\bar y))}^*\Phi_M^{-1}\otimes\Tan_y\Phi&
\\
\Tan_{\bar y}^*M\otimes{\rm V}_y(\pi)& $\rightarrowfill$ &
\Tan_{\bar\pi^1 (j^1\Phi (\bar y))}^*M\otimes{\rm V}_{\Phi (y)}(\pi)
\end{array}
$$
with the canonical identifications.
\end{itemize}

\subsubsection{Relation with the geometric elements}

Next we analyze the relation between the canonical prolongation
of diffeomorphisms and the
contact module, the structure canonical form and the canonical
isomorphisms, proving that all of them are invariant
under the action of these prolongations.

\begin{prop}
Let $\Phi\colon E\to E$ be a diffeomorphism and
$j^1\Phi\colon J^1E\to J^1E$ its canonical prolongation. Then
$$
1. \quad(j^1\Phi )^*{\cal M}_c={\cal M}_c
\hspace{2cm}
2. \quad (j^1\Phi )^*\theta = \theta
$$
\end{prop}
( {\sl Proof} )\quad
\ben
\item
Consider $\alpha\in {\cal M}_c$. In order to see that
$(j^1\Phi )^*\alpha\in {\cal M}_c$ we have to prove that it vanishes
on every holonomic section. Let $\phi\colon U \to E$ be
a local section and $j^1\phi\colon U \to J^1E$
its canonical prolongation, we have
$$
(j^1\phi)^*((j^1\Phi )^*\alpha)=
(j^1\Phi\circ j^1\phi )^*\alpha=
(j^1(\Phi\circ \phi\circ\Phi_M^{-1})\circ\Phi_M)^*\alpha=
(\Phi_M^*\circ j^1(\Phi\circ \phi\circ\Phi_M^{-1})^*)\alpha=0
$$
because $j^1(\Phi\circ \phi\circ\Phi_M^{-1})^*\alpha=0$
since $\alpha\in {\cal M}_c$ and $\Phi\circ \phi\circ\Phi_M^{-1}$
is a section of $\pi$.

We have proved that $(j^1\Phi)^*{\cal M}_c\subset{\cal M}_c$
and the converse inclusion is a direct consequence
of the third item of proposition \ref{ppc}.
\item
Consider $\theta$ as a $\Cinfty (J^1E)$-linear map
between $\vf (J^1E)$ and $\Gamma (J^1E,\pi^{1^*}{\rm V}(\pi))$.
Then we have
$$
\begin{array}{ccc}
\vf (J^1E) & $\rightarrowfill$ & \vf (J^1E)
\\
& (j^1\Phi )_* &
\\
(j^1\Phi )^*\theta \ \Big\downarrow & &\Big\downarrow \ \theta
\\
& (j^1\Phi )_* &
\\
\Gamma (J^1E,\pi^{1^*}{\rm V}(\pi)) & $\rightarrowfill$ &
\Gamma (J^1E,\pi^{1^*}{\rm V}(\pi))
\end{array}
$$
therefore
$$
(j^1\Phi )^*\theta = (j^1\Phi )^{-1}_*\circ\theta\circ (j^1\Phi)_*
$$
and  $(j^1\Phi )^*\theta =\theta$ is equivalent to
\beq
(j^1\Phi )_*\circ\theta = \theta\circ (j^1\Phi)_*
\label{aux1}
\eeq

In order to prove this, let $\bar y\in J^1E$,
$u\in\Tan_{\bar y}J^1E$ and $\phi\in\bar y$, ($\phi\colon U\to E$).
Taking into account the action of $j^1\Phi$ on $\pi^{1^*}{\rm V}(\pi)$,
we have
\beann
(\Tan j^1\Phi\circ\theta )(\bar y;u) &=&
\Tan j^1\Phi (\theta (\bar y;u))=
\Tan j^1\Phi (\d^v_y\phi(\Tan_{\bar y}\pi^1 (u)))
\\ &=&
(j^1\Phi)_*(\Tan_{\bar y}\pi^1 (u)-\Tan_{\bar y}(\phi\circ\bar\pi^1 )(u))
\\ &=&
\Tan_{\bar y}(\Phi\circ\pi^1)(u)-
\Tan_{\bar y}(\Phi\circ \phi\circ\bar\pi^1 )(u)
\eeann
however, since a representative of $(j^1\Phi)(\bar y)$ is
$\psi=\Phi\circ \phi\circ\Phi^{-1}_M$, we obtain
\beann
(\theta\circ \Tan j^1\Phi )(\bar y;u) &=&
\theta ((j^1\Phi )(\bar y);\Tan_{\bar y}(u))
\\ &=&
(\d^v_{\Phi (y)}\psi)((\Tan_{(j^1\Phi )(\bar y)}\pi^1\circ
\Tan_{\bar y}j^1\Phi )(u))
\\ &=&
(\d^v_{\Phi (y)}\psi)(\Tan_{\bar y}(\pi^1\circ j^1\Phi )(u))
\\ &=&
\Tan_{\bar y}(\pi^1\circ j^1\Phi )(u)-(\Tan_y(\psi\circ\pi)\circ
\Tan_{\bar y}(\pi^1\circ j^1\Phi ))(u)
\\ &=&
\Tan_{\bar y}(\pi^1\circ j^1\Phi )(u)-
\Tan_y(\Phi\circ \phi\circ\Phi^{-1}_M\circ\bar\pi^1\circ j^1\Phi )(u)
\\ &=&
\Tan_{\bar y}(\pi^1\circ j^1\Phi )(u)-\Tan_y(\Phi\circ \phi\circ\bar\pi^1 )(u)
\eeann
\een
\qed

Conversely, we can characterize the diffeomorphisms in $J^1E$
which are canonical prolongations as follows:

\begin{prop}
Let $\Phi\colon E\to E$ be a diffeomorphism of
$\pi$-fiber bundles and $\Phi_M\colon M\to M$ the diffeomorphism
induced on the basis.
Let $\tilde\Phi\colon J^1E\to J^1E$ be a diffeomorphism
verifying the following conditions:
$$
1. \quad \pi^1\circ\tilde\Phi =\Phi\circ\pi^1
\hspace{2cm}
2. \quad \tilde\Phi^*\theta =\theta
$$
Then $\tilde\Phi =j^1\Phi$.
\end{prop}
( {\sl Proof} )\quad
Let $\bar y\in J^1E$ and $\phi\in\bar y$, with
$\phi\colon U\to E$, $U\subset M$. The second condition implies that
$(j^1\phi )^*(\tilde\Phi^*\theta )=(j^1\phi )^*\theta =0$,
then
$$
({\Phi_M^{-1}}^*\circ (j^1\phi )^*\circ\tilde\Phi^*)\theta =
(\tilde\Phi\circ j^1\phi\circ\Phi_M^{-1})^*\theta =0
$$
However,
\beann
\bar\pi^1\circ\tilde\Phi\circ j^1\phi\circ\Phi_M^{-1}
&=&
\pi\circ\pi^1\circ\tilde\Phi\circ j^1\phi\circ\Phi_M^{-1}=
\pi\circ\Phi\circ\pi^1\circ j^1\phi\circ\Phi_M^{-1}
\\ &=&
\Phi_M\circ\pi\circ\phi\circ\Phi_M^{-1}=
\Phi_M\circ{\rm Id}_U\circ\Phi_M^{-1}={\rm Id}_{\Phi_M(U)}
\eeann
then the map
$\tilde\Phi\circ j^1\phi\circ\Phi_M^{-1}\colon\Phi_M(U)\to J^1E$
is a section of $\bar\pi^1$ such that $\theta$ vanishes on it,
so it is a holonomic section and there exists
$\varphi\colon\Phi_M(U)\to E$,
with $\tilde\Phi\circ j^1\phi\circ\Phi_M^{-1}=j^1\varphi$.
Therefore
$$
\varphi =
\pi^1\circ j^1\varphi =\pi^1\circ\tilde\Phi\circ j^1\phi\circ\Phi_M^{-1}=
\Phi\circ\pi^1\circ j^1\phi\circ\Phi_M^{-1}=\Phi\circ\phi\circ\Phi_M^{-1}
$$
and we have
\beann
\tilde\Phi (\bar y)&=&
\tilde\Phi (j^1\phi (\bar\pi^1 (\bar y)))=
(\tilde\Phi\circ j^1\phi\circ\Phi_M^{-1})(\Phi_M(\bar\pi^1 (\bar y)))
\\ &=&
(j^1\varphi )(\Phi_M(\bar\pi^1 (\bar y)))=
j^1(\Phi\circ j^1\phi\circ\Phi_M^{-1})(\Phi_M(\bar\pi^1 (\bar y)))=
(j^1\Phi )(\bar y)
\eeann
so $\tilde\Phi =j^1\Phi$.
\qed

The last two results can be summarized in the following assertion:

\begin{prop}
Let $\tilde\Phi\colon J^1E\to J^1E$ be a diffeomorphism of
fiber bundles over $E$ and over $M$,
which induces $\Phi\colon E\to E$ and
$\Phi_M\colon M\to M$.
The necessary and sufficient condition for $\tilde\Phi$
to be the canonical prolongation to $J^1E$
of the induced diffeomorphism $\Phi$ on $E$ is that
$\tilde\Phi^*\theta =\theta$.
\end{prop}

Finally, for the vertical endomorphisms we have:

\begin{prop}
Let $\Phi\colon E\to E$ be a diffeomorphism and
$j^1\Phi\colon J^1E\to J^1E$ its canonical prolongation. Then
$$
1. \quad (j^1\Phi )^*{\cal S} = {\cal S}
\hspace{2cm}
2. \quad (j^1\Phi )^*{\cal V} = {\cal V}
$$
\end{prop}
( {\sl Proof} )\quad
We are going to interpret $(j^1\Phi )^*{\cal S}$ as a morphism of sections
over $J^1E$. Remember that, if $\bar y\in J^1E$ with $\bar y\mapsto y\mapsto x$
and $\alpha\otimes v\in{\rm V}_y(\pi )\otimes\Tan_x^*M$, then
${\cal S}(\bar y;\alpha\otimes v )=D_{\alpha\otimes v}(\bar y)$,
the directional derivative in $J^1_yE$
at the point $\bar y$ in the direction given by $\alpha\otimes v$.
According to this we can display the following diagram
$$
\begin{array}{ccc}
\Gamma (J^1E,\bar\pi^{1^*}\Tan^*M)\otimes\Gamma (J^1E,\pi^{1^*}{\rm V}(\pi)) &
$\rightarrowfill$
& \Gamma (J^1E,\bar\pi^{1^*}\Tan^*M)\otimes\Gamma (J^1E,\pi^{1^*}{\rm V}(\pi))
\\
& (j^1\Phi^{-1})^*\otimes (j^1\Phi )_* &
\\
(j^1\Phi )^*{\cal S}  \Big\downarrow & &\Big\downarrow  {\cal S}
\\
& (j^1\Phi )_* &
\\
\Gamma (J^1E,{\rm V}(\pi^1 )) & $\rightarrowfill$ &
\Gamma (J^1E,{\rm V}(\pi^1 ))
\end{array}
$$
therefore, as in the above proposition of invariance of the structure form
$\theta$, we have
$$
(j^1\Phi )^*{\cal S}=
((j^1\Phi^{-1})^*\otimes(j^1\Phi)_*)^{-1}\circ {\cal S}\circ (j^1\Phi )^*
$$
and the statement of the theorem is equivalent to
$$
{\cal S}\circ ((j^1\Phi^{-1})^*\otimes(j^1\Phi)_*)=(j^1\Phi )_*\circ {\cal S}
$$

We have
\beann
({\cal S}\circ ((j^1\Phi^{-1})^*\otimes(j^1\Phi)_*))(\bar y;\alpha\otimes v )
&=&
{\cal S}((j^1\Phi )(\bar y);\Tan^*j^1\Phi^{-1}\alpha\otimes \Tan j^1\Phi v)
\\ &=&
D_{\Tan^*j^1\Phi^{-1}\alpha\otimes\Tan j^1\Phi v}(j^1\Phi (\bar y))
\eeann
However
$$
((j^1\Phi )_*\circ {\cal S})(\bar y;\alpha\otimes v )=
(\Tan_{\bar y}j^1\Phi )(D_{\alpha\otimes v}(\bar y))
$$
but, if $f\in\Cinfty (J^1E)$, then
$$
D_{\Tan^*j^1\Phi^{-1}\alpha\otimes\Tan j^1\Phi v}(j^1\Phi (\bar y))f=
\lim_{t\to 0}\frac{1}{t}(f(j^1\Phi )(\bar y)+
t(\Tan^*j^1\Phi^{-1}\alpha\otimes\Tan j^1\Phi v))-
f(j^1\Phi )(\bar y)))
$$
and
$$
(\Tan j^1\Phi (D_{\alpha\otimes v}(\bar y))f=
\lim_{t\to 0}\frac{1}{t}(f((j^1\Phi )(\bar y+t\alpha\otimes v))-
f(j^1\Phi )(\bar y)))
$$
then it must be proved that
$$
(j^1\Phi )(\bar y+t\alpha\otimes v)=(j^1\Phi )(\bar y)
+t(\Tan^*j^1\Phi^{-1}\alpha\otimes\Tan j^1\Phi v)
$$
But this is proposition \ref{canprol}.

Since ${\cal V}=\inn ({\cal S})\theta$,
the second item obviously holds.
\qed

\subsubsection{Prolongation of projectable vector fields}

Let $Z\in\vf (E)$ a $\pi$-projectable vector field
and $\tau_t\colon W\to W$, ($W\subset E$ being an open set),
a local one-parameter group of diffeomorphisms of $Z$.
$\tau_t$ are bundle-diffeomorphisms, hence there exist
their prolongations $j^1\tau_t$ to $J^1E$. From the properties of $\tau_t$,
we deduce that $j^1\tau_t$ is also a local one-parameter group of
diffeomorphisms defined in $(\pi^1)^{-1}(W)$, so they are associated
to a vector field. Since the prolongations can be done on every open set
where a local one-parameter group of diffeomorphisms of $Z$ exists,
we can extend this vector field everywhere in its domain.

\begin{definition}
Let $Z\in\vf (E)$ be a $\pi$-projectable vector field.
The {\rm canonical prolongation} of $Z$
is the vector field $j^1Z\in\vf(J^1E)$ whose associated
local one-parameter groups of diffeomorphisms are
the canonical prolongations $j^1\tau_t$ of the local
one-parameter groups of diffeomorphisms of $Z$.
\end{definition}

The most relevant properties of these prolongations are:

\begin{prop}
Let $Z\in\vf (E)$ be a $\pi$-projectable vector field.
Then:
\ben
\item
$j^1Z$ is well defined and it is a $\pi^1$-projectable
vector field, satisfying that $\pi^1_*(j^1Z)=Z$.
\item
$\Lie (j^1Z){\cal M}_c\subset {\cal M}_c$.
\item
The canonical prolongation of $Z$
is the only vector field $j^1Z\in\vf(J^1E)$
verifying the two following conditions:
\ben
\item
$j^1Z$ is $\pi^1$-projectable and $\pi^1_*(j^1Z)=Z$.
\item
$j^1Z$ lets the contact module invariant:
$\Lie (j^1Z){\cal M}_c\subset {\cal M}_c$.
\een
\een
\label{provf}
\end{prop}
\proof
\ben
\item
It is a direct consequence of the theorem of existence and unicity
of the local one-parameter groups.
\item
Since the canonical prolongations of diffeomorphisms preserve the
contact module, the result follows straightforwardly.
\item
It suffices to see that these conditions allow us to calculate
the coordinates of $j^1Z$ in any local canonical system of $J^1E$.
So, let $(x^\mu ,y^A,v^A_\mu)$ be a local natural system in
$(\pi^1)^{-1}(W)\subset J^1E$, and
\dst Z=\alpha^\mu\derpar{}{x^\mu}+\beta^A\derpar{}{y^A}\)
a $\pi$-projectable vector field; that is,
\dst \derpar{\alpha^\mu}{y^A}=0\) .
Since $j^1Z$ is $\pi^1$-projectable and $\pi^1_*(j^1Z)=Z$, we have
$$
j^1Z=\alpha^\mu\derpar{}{x^\mu}+\beta^A\derpar{}{y^A}
+\gamma^A_\mu\derpar{}{v^A_\mu}
$$
and the only problem is to calculate the coefficients $\gamma^A_\mu$.
But we know that $j^1Z$ preserves the contact module and this is
generated by the forms $\theta^A=\d y^A-v^A_\mu\d x^\mu$,
then we have
$$
\Lie (j^1Z)\theta^A=\zeta^A_B\theta^B
$$
where $\zeta^A_B\in\Cinfty (W)$. We obtain
\beann
\zeta^A_B\theta^B&=&\Lie (j^1Z)\theta^A=
(\d\inn (j^1Z)+\inn(j^1Z)\d )\theta^A
\\ &=&
\d (\beta^A-v^A_\mu\alpha^\mu )+\inn (j^1Z)(\d v^A_\mu\wedge\d x^\mu )
\\ &=&
\derpar{\beta^A}{x^\mu}\d x^\mu+\derpar{\beta^A}{y^B}\d y^B
-\alpha^\mu\d v^A_\mu-v^A_\mu\derpar{\alpha^\mu}{x^\rho}\d x^\rho
-\gamma^A_\mu\d x^\mu+\alpha^\mu\d v^A_\mu
\\ &=&
\left(\derpar{\beta^A}{x^\mu}-v^A_\rho\derpar{\alpha^\rho}{x^\mu}
-\gamma^A_\mu\right)\d x^\mu+\derpar{\beta^A}{y^B}\d y^B
\eeann
Then, by identification of terms
$$
\zeta^A_B=\derpar{\beta^A}{y^B} \quad , \quad
-\zeta^A_Bv^B_\mu=\derpar{\beta^A}{x^\mu}-v^A_\rho\derpar{\alpha^\rho}{x^\mu}
-\gamma^A_\mu
$$
therefore
$$
\gamma^A_\mu=
\derpar{\beta^A}{x^\mu}-v^A_\rho\derpar{\alpha^\rho}{x^\mu}+
v^B_\mu\derpar{\beta^A}{y^B}
$$
which completes the local coordinate expression of $j^1Z$.
\een
\qed

As a particular case,
if $Z\in\vf (E)$ is a vertical vector field;
that is, ${\pi}_*Z=0$, then if $\tau_t$ is a local one-parameter group
associated with $Z$, it induces the identity on $M$. In this case,
in a local natural system of coordinates,
\dst Z=\beta^A\derpar{}{y^A}\) and its canonical prolongation is
$$
j^1Z=\beta^A\derpar{}{y^A}+
\left(\derpar{\beta^A}{x^\mu}+\derpar{\beta^A}{y^B}v^B_\mu\right)
\derpar{}{v^A_\mu}
$$

\subsubsection{Generalization to non-projectable vector fields}

If $Z\in\vf (E)$ is not a $\pi$-projectable vector field,
we can define its prolongation by means of the characterization
given in proposition \ref{provf}. Thus we define:

\begin{definition}
Let $Z\in\vf (E)$ be a vector field.
The {\rm canonical prolongation} of $Z$
is the only vector field $j^1Z\in\vf(J^1E)$ verifying the following conditions:
\ben
\item
$j^1Z$ is $\pi^1$-projectable and $\pi^1_*(j^1Z)=Z$.
\item
The contact module is invariant under the action of $j^1Z$:
$\Lie (j^1Z){\cal M}_c\subset {\cal M}_c$.
\een
A vector field $X\in\vf (J^1E)$ satisfying these conditions is called an
{\rm infinitesimal contact transformation}.
\end{definition}

In a local system of coordinates, $(x^\mu ,y^A,v^A_\mu)$,
with \dst Z=\alpha^\mu\derpar{}{x^\mu}+\beta^A\derpar{}{y^A}\) ,
doing an analogous calculation as in the above section
and taking into account that, in general,
\dst \derpar{\alpha^\mu}{y^A}\not=0\) ,
we obtain finally
\beann
j^1Z &=& \alpha^\mu\derpar{}{x^\mu}+\beta^A\derpar{}{y^A}+
\left(\derpar{\beta^A}{x^\mu}-
v^A_\rho\left(\derpar{\alpha^\rho}{x^\mu}+
v^B_\mu\derpar{\alpha^\rho}{y^B}\right)+v^B_\mu\derpar{\beta^A}{y^B}\right)
\derpar{}{v^A_\mu}
\eeann

For further points of view on this subject,
see \cite{Gc-74} or \cite{GIMMSY-mm}.

\section{Lagrangian formalism of first-order field theories}

\subsection{Lagrangian systems. Variational problems and critical sections}

This section is devoted to the introduction of {\sl lagrangian densities}
and, using them, to defining the geometrical objects that allow us
to study the dynamical behaviour of field theories.

Let us consider again the initial situation
$\pi\colon E\to M$ and the jet bundle $\pi^1\colon J^1E\to E$,
$M$ being an orientable manifold.

\subsubsection{Lagrangian densities}

Dynamics of classical field theories is given
from lagrangian densities. A {\sl lagrangian density}
depends on the fields and their first derivatives
with respect to the space-time variables, and it can be integrated
on the images of the fields.

First of all, it is important to point out that
the following $\Cinfty (J^1E)$-modules are canonically isomorphic:
\ben
\item
$\Cinfty (J^1E)\otimes_M\df^{n+1}(M)$.
\item
$\Gamma (J^1E,\bar\pi^{1^*}\Lambda^{n+1}\Tan^*M)$.
\item
The set of semibasic $(n+1)$-forms on $J^1E$ with respect to the projection
$\bar\pi^1$:
$$
\{ \alpha\in\df^{n+1}(J^1E)\ |\ \inn (X)\alpha =0\ ,\
\forall X\in\Gamma (J^1E,{\rm V}(\bar\pi^1))\}
$$
\item
The sections of the bundle $\Lambda^{n+1}\Tan^*M$
along the mapping $\bar\pi^1\colon J^1E\to M$.
\een
(Proofs of this proposition can be found in
\cite{GHV-72} and \cite{LM-87}).
\qed

This assertion allow us a free ride on the notion of
{\sl lagrangian density}. Then we define:

\begin{definition}
A {\rm lagrangian density} $\Lag$ is a $\bar\pi^1$-semibasic
$(n+1)$-form in $J^1E$.
\end{definition}

{\bf Remarks}:
\ben
\item
After the above proposition, the lagrangian densities will be considered
indistinctly as elements of any one of those $\Cinfty (J^1E)$-modules
when necessary.
\item
The expression of a lagrangian density $\Lag$ is
$$
\Lag = f_i\otimes\bar\pi^{1^*}\omega^i =f_i\otimes\omega^i
\quad ; \quad
f_i\in\Cinfty (J^1E)\ ,\ \omega^i\in\df^{n+1}(M)
$$
Since we suppose that $M$ is an oriented manifold,
suppose $\omega\in\df^{n+1}(M)$ is a fixed volume $(n+1)$-form on $M$,
then from now on we will write a lagrangian density as
$$
\Lag =\lag\omega \quad ; \quad
\lag\in\Cinfty (J^1E)\ ,\ \omega\in\df^{n+1}(M)
$$
If $(x^\mu ,y^A,v^A_\mu )$ is a natural system of coordinates,
we use the expression
$$
\Lag = \lag (x^\mu ,y^A,v^A_\mu )\d x^1\wedge\ldots\wedge\d x^n\wedge\d x^0
\equiv\lag\d^{n+1}x
$$
Then:

\begin{definition}
The function $\lag$ is called the {\rm lagrangian function}
associated with $\Lag$ and $\omega$.
\end{definition}
\item
Bearing this in mind, the exterior differential of $\Lag$
is $\d\Lag =\d\lag\wedge\omega$, which is an element of
$\df^{n+2}(J^1E)$ or an element of
$\df^1(J^1E)\wedge\Gamma (J^1E,\bar\pi^{1^*}\Lambda^{n+1}\Tan^*M)$,
with the exterior product as $\Cinfty (J^1E)$-modules.

Observe that, if $X_1,X_2$ are $\pi^1$-vertical vector fields,
then $\inn (X_1)\inn(X_2)\d\Lag =0$.
\een

Finally, we establish the following terminology:

\begin{definition}
A {\rm lagrangian system} is a pair $\ls$ where
\ben
\item
$\pi\colon E\to M$ is a differentiable bundle with $M$ an orientable
manifold.
\item
$\Lag$ is a lagrangian density.
\een
\end{definition}

\subsubsection{Lagrangian canonical forms. Properties}

Now we can define the canonical forms induced by
a lagrangian density. Given a lagrangian density
$\Lag$ which we write in the form $\Lag =\lag\omega$, we are going to associate
it with some differential forms in order to use them
for studying variational problems.

The standpoint objects are
\beann
{\cal V}&\in&\df^1(J^1E)\otimes\Gamma (J^1E,{\rm V}(\pi^1 ))\otimes
\Gamma (J^1E,\bar\pi^{1^*}\Tan M)
\\
\d\Lag&\in&\df^1(J^1E)\wedge\Gamma (J^1E,\bar\pi^{1^*}\Lambda^{n+1}\Tan^*M)
\eeann
and we can contract the factors
$\Gamma (J^1E,{\rm V}(\pi^1 ))$ of ${\cal V}$
with $\df^1(J^1E)$ of $\d\Lag$, and
$\Gamma (J^1E,\bar\pi^{1^*}\Tan M)$ with
$\Gamma (J^1E,\bar\pi^{1^*}\Lambda^{n+1}\Tan^*M)$,
obtaining a factor in $\Gamma (J^1E,\bar\pi^{1^*}\Lambda^n\Tan^*M)$.
So we have an element of
$\df^1(J^1E)\wedge\Gamma (J^1E,\bar\pi^{1^*}\Lambda^n\Tan^*M)$
which is also an element of $\df^{n+1}(J^1E)$
in a natural way. Then

\begin{definition}
Let $\ls$ be a lagrangian system.
\ben
\item
The $(n+1)$-form in $J^1E$
$$
\vartheta_{\Lag}:=\inn ({\cal V})\d\Lag
$$
(obtained after doing the above-mentioned contractions) is called the
{\rm lagrangian canonical form} associated with the lagrangian density $\Lag$.
\item
The $(n+1)$-form in $J^1E$
$$
\Theta_{\Lag}:=\vartheta_{\Lag} +\Lag
$$
is called the
{\rm Poincar\'e-Cartan $(n+1)$-form} associated with the lagrangian density
$\Lag$.
\item
The $(n+2)$-form in $J^1E$
$$
\Omega_{\Lag}:= -\d\Theta_{\Lag}
$$
is called the {\rm Poincar\'e-Cartan $(n+2)$-form}
associated with the lagrangian density $\Lag$.
\een
\end{definition}

In a natural system of coordinates $(x^\mu ,y^A,v^A_\mu )$
the expressions of these elements are:
\beann
{\cal V}&=&\theta^A\otimes\derpar{}{v^A_\mu}\otimes\derpar{}{x^\mu}
\\
\theta^A&=&\d y^A-v^A_\rho\d x^\rho
\\
\d\Lag &=&\d\lag\wedge\omega
\\
\omega &=&\d x^1\wedge\ldots\d x^n\wedge\d x^0\equiv\d^{n+1}x
\\
\vartheta_{\Lag}:=\inn ({\cal V})\d\Lag&=&
\derpar{\lag}{v^A_\mu}\theta^A\wedge\inn\left(\derpar{}{x^\mu}\right)\omega=
\derpar{\lag}{v^A_\mu}\theta^A\wedge\d^nx_\mu
\\ &\equiv&
\derpar{\lag}{v^A_\mu}\theta^A\wedge(-1)^{\mu +1}\d x^1\wedge\ldots
\d x^{\mu -1}\wedge\d x^{\mu +1}\wedge\ldots\d x^0
\\ &=&
\derpar{\lag}{v^A_\mu}\d y^A\wedge\d^nx_\mu-
\derpar{\lag}{v^A_\mu}v^A_\mu\d^{n+1}x
\\
\Theta_{\Lag}&=&\derpar{\lag}{v^A_\mu}\d y^A\wedge\d^nx_\mu -
\left(\derpar{\lag}{v^A_\mu}v^A_\mu -\lag\right)\d^{n+1}x
\\
\Omega_{\Lag}&=&
-\d\left(\derpar{\lag}{v^A_\mu}\right)\wedge\d y^A\wedge\d^nx_\mu +
\d\left(\derpar{\lag}{v^A_\mu}v^A_\mu -\lag\right)\wedge\d^{n+1}x
\eeann

As a direct consequence of the definition and the above local
expressions, the relation between these canonical forms and
the prolongation of sections is the following:

\begin{prop}
Let $\phi\colon U\subset M\to E$ be a local section of $\pi$.
If $\Lag$ is a lagrangian density then
\ben
\item
$(j^1\phi )^*\vartheta_{\Lag} =0$.
\item
$(j^1\phi )^*\Theta_{\Lag} =(j^1\phi )^*\Lag$.
\een
\end{prop}

\subsubsection{Density of lagrangian energy}

In the last paragraph we have seen how,
in the local expression of the Poincar\'e-Cartan $(n+2)$-form,
appears the factor
\beq
{\rm E}_{\lag}\equiv\derpar{\lag}{v^A_\mu}v^A_\mu -\lag
\label{clen}
\eeq
which is recognized as the classical expression of the
{\sl lagrangian energy} associated with the lagrangian function $\lag$.
Is there any global function on $J^1E$ having this local expression?.
We are going to prove that, in order to give an intrinsic
construction of this function and, by extension,
of a {\sl density of lagrangian energy}, it is necessary
to choose a connection in the bundle $\pi\colon E\to M$
\footnote{
In a forthcomming work \cite{EMR-ghft},
we will prove that this density of lagrangian energy allows us
to construct a {\sl hamiltonian density} and a formulation
of the theory in hamiltonian terms.
}.
This dependence is hidden in the given local expression
because any local systems of coordinates in this bundle
induces a local connection with vanishing Christoffel symbols.

First of all, we refer to the appendix \ref{conec} in order to
take into account the basic concepts about connections
in the bundle $\pi\colon E\to M$.
Then we will see how the identifications made in this appendix
influence the vertical endomorphisms and the lagrangian forms.

Consider the bundle sequence
$$
J^1E\longrightarrow E\longrightarrow M
$$
and the bundle $\pi^{1^*}\Tan E\to J^1E$.
If $\nabla$ is a connection in $\pi\colon E\to M$
inducing the splitting
$\Tan E={\rm V}(\pi )\oplus{\rm H}(\nabla )$,
we also have that
$$
\pi^{1^*}\Tan E=
\pi^{1^*}{\rm V}(\pi )\oplus\pi^{1^*}{\rm H}(\nabla )
$$
and therefore we obtain
$$
\pi^{1^*}\Tan^*E=
\pi^{1^*}{\rm V}^*(\pi )\oplus\pi^{1^*}{\rm H}^*(\nabla )
$$
as a splitting in the dual bundle. Hence
$\pi^{1^*}{\rm V}^*(\pi )$ is identified as a subbundle of
$\pi^{1^*}\Tan^*E$.

If $\xi$ is a section of $\pi^{1^*}\Tan^*E\to J^1E$; that is,
$\xi\in\Gamma (J^1E,\pi^{1^*}\Tan^*E)$,
then $\xi$ is an element of $\df^1(J^1E)$,
that is, a $1$-form on $J^1E$, in the following way:
let $\bar y\mapsto y\mapsto x$ and $X\in{\cal X}(J^1E)$.
We write
$$
\xi (\bar y;X)=\xi_y((\Tan_{\bar y}\pi^1)(X_{\bar y}))
$$
where $\xi_y$ is the value of $\xi$ at $\bar y$, translated to
$y=\pi^1(\bar y)$ by means of the natural identification between the
corresponding fibers of
$\pi^{1^*}\Tan^*E$ and $\Tan^*E$.

In the same way, as every connection induces an injection
$\pi^{1^*}{\rm V}^*(\pi )\hookrightarrow \pi^{1^*}\Tan^*E$,
then every $\zeta\in\Gamma (J^1E,\pi^{1^*}{\rm V}^*(\pi ))$
is a $1$-form in $J^1E$. We then have
$$
\Gamma (J^1E,\pi^{1^*}{\rm V}^*(\pi ))\hookrightarrow\df^1(J^1E)
$$

Consider now the vertical endomorphisms ${\cal S}$ and ${\cal V}$,
which, as we know, are sections of the bundles
\beann
{\cal S}&\in&\Gamma (J^1E,(\pi^{1^*}{\rm V}(\pi))^*)\otimes
\Gamma (J^1E,{\rm V}(\pi^1 ))\otimes\Gamma(J^1E,\bar\pi^{1^*}\Tan M)
\\
{\cal V}&\in&\df^1(J^1E)\otimes\Gamma (J^1E,{\rm V}(\pi^1 ))
\otimes\Gamma (J^1E,\bar\pi^{1^*}\Tan M)
\eeann
(remember that
$(\pi^{1^*}{\rm V}(\pi))^*=\pi^{1^*}{\rm V}^*(\pi)$).
Taking into account the above identifications,
the following operation holds
$$
{\cal S}-{\cal V}\in\df^1(J^1E)\otimes\Gamma (J^1E,{\rm V}(\pi^1 ))
\otimes\Gamma (J^1E,\bar\pi^{1^*}\Tan M)
$$
Consequently we are able to achieve the contraction of this difference with
differential forms in $J^1E$ such as $\d\Lag$,
obtaining
$$
\inn ({\cal S}-{\cal V})\d\Lag\in\df^{n+1}(J^1E)
$$
which is a semibasic form with respect to $\pi^1$. Then we can define:

\begin{definition}
Let $\ls$ be a lagrangian system and
$\nabla$ be a connection in $\pi\colon E\to M$. The
{\rm density of lagrangian energy} associated with the lagrangian density
$\Lag$ and the connection $\nabla$ is the semibasic $(n+1)$-form
in $J^1E$ given by
$$
{\cal E}^\nabla_{\Lag}:=\inn ({\cal S}-{\cal V})\d\Lag -\Lag
$$
\end{definition}

As in the case of the lagrangian density,
the expression of the density of lagrangian energy is
$$
\del = F_i\otimes\bar\pi^{1^*}\omega^i =F_i\otimes\omega^i
\quad ; \quad
F_i\in\Cinfty (J^1E)\ ,\ \omega^i\in\df^{n+1}(M)
$$
and, if $\omega\in\df^{n+1}(M)$ is the volume $(n+1)$-form on $M$,
from now on we will write the density of lagrangian energy as
$$
\del =\fel\bar\pi^{1^*}\omega\equiv\fel\omega \quad ; \quad
\fel\in\Cinfty (J^1E)\ ,\ \omega\in\df^{n+1}(M)
$$
Then:

\begin{definition}
The function $\fel$ is called the {\rm lagrangian energy function}
associated with $\Lag$, $\nabla$ and $\omega$.
\end{definition}

In order to see the local expressions of these elements,
let $(x^\mu ,y^A,v^A_\mu )$ be a natural system of coordinates.
We have
\beann
\nabla&=&
\d x^\mu\otimes\left(\derpar{}{x^\mu}+{\mit\Gamma}_\mu^A\derpar{}{y^A}\right)
\\
\del &=&\fel\d^{n+1}x
\\
{\cal S}&=&(\d y^A-{\mit\Gamma}^A_\mu\d x^\mu )
\otimes\derpar{}{v^A_\nu}\otimes\derpar{}{x^\nu}
\\
{\cal V}&=&\left(\d y^A-v^A_\mu\d x^\mu\right)\otimes
\derpar{}{v^A_\nu}\otimes\derpar{}{x^\nu}
\eeann
where, in the expression of ${\cal S}$, we have taken
$\{\d y^A-{\mit\Gamma}^A_\mu\d x^\mu\}$ as a local basis of
the sections of $\pi^{1^*}{\rm V}^*(\pi )$ over $J^1E$,
associated with the connection $\nabla$.
This basis is well defined because $\pi^{1^*}{\rm V}^*(\pi )$
is the incident submodule to $\pi^{1^*}{\rm H}(\nabla )$,
which is locally generated by
\dst\left\{\derpar{}{x^\mu}+{\mit\Gamma}^A_\mu\derpar{}{y^A}\right\}\) .
Hence
$$
{\cal S}-{\cal V}=(v^A_\mu-{\mit\Gamma}^A_\mu )\d x^\mu
\otimes\derpar{}{v^A_\nu}\otimes\derpar{}{x^\nu}
$$
and then the density of lagrangian energy is
$$
{\cal E}_{\Lag}^\nabla:=\inn ({\cal S}-{\cal V})\d\Lag -\Lag =
\left(\derpar{\lag}{v^A_\mu} (v^A_\mu-{\mit\Gamma}^A_\mu )-\lag\right)
\d^{n+1}x
$$
where
$$
{\rm E}_{\Lag}^\nabla=
\derpar{\lag}{v^A_\mu} (v^A_\mu-{\mit\Gamma}^A_\mu )-\lag
$$
is the expression of the lagrangian energy.
As one may observe, the definition of the density of lagrangian energy
associated with $\Lag$ and $\nabla$ depends only on these elements.

{\bf Remark}:
\bit
\item
If the bundle $\pi\colon E\to M$ is trivial,
that is $E=F\times M$, then we have a natural connection on $E$
given by the splitting $\Tan E=\Tan F\times\Tan M$,
and the density of lagrangian energy associated with this
natural connection is just (\ref{clen}),
since ${\mit\Gamma}^A_\mu =0$.
\eit

As a final remark, we will compare the densities of lagrangian energy
corresponding to two different connections.

First of all, as is described in appendix \ref{conec},
given a connection $\nabla$
and $\gamma\in\Gamma (E,\pi^*\Tan^*M)\otimes\Gamma (E,{\rm V}(\pi ))$,
then $\nabla+\gamma$ is also a connection and, conversely,
if $\nabla$ and $\nabla'$ are connections then
$\nabla-\nabla'=\gamma\in
\Gamma (E,\pi^*\Tan^*M)\otimes\Gamma (E,{\rm V}(\pi ))$;
that is, the difference between two connections is not a connection.
Hence, for a given connection $\nabla$,
the construction of the density of lagrangian energy allows us to define a map
\beann
\Gamma (E,\pi^*\Tan^*M)\otimes\Gamma (E,{\rm V}(\pi )) &
\longrightarrow &\Gamma (J^1E,\bar\pi^{1^*}\Lambda^{n+1}\Tan^*M)
\\
\gamma & \mapsto & {\cal E}^{\nabla+\gamma}_{\Lag}-{\cal E}^{\nabla}_{\Lag}
\eeann
which is a $\Cinfty (E)$-morphism of modules. This is the so-called
{\sl Legendre transformation form} in \cite{Gc-74}.

In another way, if we denote by
${\cal S}^\nabla$ and ${\cal S}^{\nabla +\gamma}$
the ``vertical endomorphisms corresponding to $\nabla$ and $\nabla +\gamma$''
respectivelly, we have
$$
{\cal E}^{\nabla+\gamma}_{\Lag}-{\cal E}^{\nabla}_{\Lag} =
\inn ({\cal S}^{\nabla +\gamma}-{\cal S}^\nabla )\d\Lag
$$

In a local system of coordinates we have
\beann
\gamma &=&\gamma^A_\mu\d x^\mu\otimes\derpar{}{y^A}
\\
{\cal E}^{\nabla+\gamma}_{\Lag}-{\cal E}^{\nabla}_{\Lag} &=&
-\gamma^A_\mu\derpar{\lag}{v^A_\nu}\d^{n+1}x
\eeann
Note that this difference does not depend on the connection
but only on $\gamma$.

This expression of the difference of densities of lagrangian energy
corresponding to two different connections is recognized by stating that
the density of lagrangian energy depends linearly on the connection.

\subsubsection{Variational problem associated with a lagrangian density}

If $\phi\colon M\to E$ is a section of
$\pi$ and $j^1\phi\colon M\to J^1E$ is its canonical prolongation,
since $\Lag$ is a $(n+1)$-form on $J^1E$, then
$(j^1\phi)^*\Lag\in\df^{n+1}(M)$.
Taking this into account we can define:

\begin{definition}
Let $\ls$ be a lagrangian system.
Let $\Gamma_c(M,E)$ be the set of compact supported sections of
$\pi$ and consider the map
$$
\begin{array}{ccccc}
{\bf L}&\colon&\Gamma_c(M,E)&\longrightarrow&\Real
\\
& &\phi&\mapsto&\int_M(j^1\phi)^*\Lag
\end{array}
$$
The {\rm variational problem} posed by the lagrangian density $\Lag$
is the problem of searching for the critical (or stationary) sections
of the functional ${\bf L}$.
\end{definition}

{\bf Remark}:
\begin{itemize}
\item
The sections must be ``stationary'' with respect to the variations of $\phi$
given by $\phi_t =\tau_t\circ \phi$, where $\{\tau_t\}$
is a local one-parameter group of any
vector field $Z\in\vf (E)$ which is $\pi$-vertical; that is,
$$
\frac{\d}{\d t}\Big\vert_{t=0}\int_M(j^1\phi_t)^*\Lag = 0
$$
$Z$ is required to be vertical in order to assure that
$\phi_t$ is a section defined in the same open set as $\phi$
(except for some problems in the domain of $\tau_t$),
since $\tau_t$ induces the identity on $M$.
\end{itemize}

\subsubsection{Characterization of critical sections}

In this paragraph we suppose that a lagrangian density
$\Lag$ is given and we are going to study the conditions for a section
of $\pi$ with compact support to be a critical section of
the variational problem associated with this lagrangian density.
For other approaches see \cite{GS-73}. \cite{Gc-74}, \cite{GIMMSY-mm}.

\begin{teor}
Let $\ls$ be a lagrangian system.
The following assertions on a section $\phi\in\Gamma_c(M,E)$
are equivalent
\footnote{
In the Lie derivatives and contractions, $\Lag$ must be understood as a
$(n+1)$-form on $J^1E$.
}.
\ben
\item
$\phi$ is a critical section for the variational problem posed by $\Lag$.
\item
\dst\int_M(j^1\phi)^*\Lie (j^1Z)\Lag = 0\) ,
for every $Z\in\vf^{{\rm V}(\pi)}(E)$.
\item
\dst\int_M(j^1\phi)^*\Lie (j^1Z)\Theta_{\Lag} = 0\) ,
for every $Z\in\vf^{{\rm V}(\pi)}(E)$.
\item
$(j^1\phi)^*\inn (j^1Z)\Omega_{\Lag} = 0$,
for every $Z\in\vf^{{\rm V}(\pi)}(E)$.
\item
$(j^1\phi)^*\inn (X)\Omega_{\Lag} = 0$,
for every $X\in\vf (J^1E)$.
\item
If $(x^\mu ,y^A,v^A_\mu )$ is a local natural system of coordinates
and $\Lag =\lag\omega =\lag\d^{n+1}x$,
the coordinates of $\phi$ in that system
satisfy the {\rm Euler-Lagrange equations}:
$$
\derpar{\lag}{y^A}\Big\vert_{j^1\phi}-
\derpar{}{x^\mu}\left(\derpar{\lag}{v_\mu^A}\Big\vert_{j^1\phi}\right) = 0
\quad ;\quad (A=1,\ldots ,N)
$$
\een
\label{equics}
\end{teor}
( {\sl Proof} )\quad
($1\ \Leftrightarrow\ 2$)\quad
If $Z\in\vf^{{\rm V}(\pi)}(E)$ and $\tau_t$
is a one-parameter local group of $Z$,
according to the results in the paragraph \ref{prol}, we have
$$
j^1(\tau_t\circ \phi)=j^1\tau_t\circ j^1\phi
$$
therefore
\beann
\frac{\d}{\d t}\Big\vert_{t=0}\int_M(j^1(\tau_t\circ \phi))^*\Lag &=&
\lim_{t\to 0}\frac{1}{t}
\left(\int_M(j^1(\tau_t\circ \phi))^*\Lag -\int_M(j^1\phi)^*\Lag\right)
\\ &=&
\lim_{t\to 0}\frac{1}{t}
\left(\int_M(j^1\phi)^*(j^1\tau_t)^*\Lag -\int_M(j^1\phi)^*\Lag\right)
\\ &=&
\lim_{t\to 0}\frac{1}{t}
\left(\int_M(j^1\phi)^*[(j^1\tau_t)^*\Lag -\Lag ]\right) =
\int_M(j^1\phi)^*\Lie (j^1Z)\Lag
\eeann
and the results follows immediatelly.

\quad\quad ($2\ \Leftrightarrow\ 3$)\quad
Taking into account that $(j^1\phi)^*{\Lag}=(j^1\phi)^*\Theta_{\Lag}$,
for every section of $\pi$, we have $\phi$ is a critical section
if, and only if, for every $Z\in\vf^{{\rm V}(\pi)}(E)$.
$$
\frac{\d}{\d t}\Big\vert_{t=0}\int_M(j^1(\tau_t\circ \phi))^*\Theta_{\Lag}
=\int_M(j^1\phi)^*\Lie (j^1Z)\Theta_{\Lag}
$$
where $\tau_t$ is the one-parameter local group associated with $Z$.

\quad\quad ($3\ \Leftrightarrow\ 4$)\quad
Taking into account the last item,
we have $\phi$ is a critical section if, and only if,
\dst\int_M(j^1\phi)^*\Lie (j^1Z)\Theta_{\Lag} =0\) ,
for every $Z\in\vf^{{\rm V}(\pi)}(E)$.
But
$$
\Lie (j^1Z)\Theta_{\Lag} =
\d\inn (j^1Z)\Theta_{\Lag}+\inn (j^1Z)\d\Theta_{\Lag} =
\d\inn (j^1Z)\Theta_{\Lag}-\inn (j^1Z)\Omega_{\Lag}
$$
and as $\phi$ has compact support, using Stoke's theorem we obtain
$$
\int_M(j^1\phi)^*\d\inn (j^1Z)\Theta_{\Lag}
=\int_M\d (j^1\phi)^*\inn (j^1Z)\Theta_{\Lag} =0
$$
and hence $\phi$ is a critical section if, and only if,
$$
\int_M(j^1\phi)^*\inn (j^1Z)\Omega_{\Lag} =0
$$
and therefore, according to the fundamental theorem of variational
calculus (see the comment at the end of the proof),
the section is critical if, and only if,
$$
(j^1\phi)^*\inn (j^1Z)\Omega_{\Lag} =0
$$

\quad\quad ($4\ \Leftrightarrow\ 5$)\quad
The sufficency is a consequence of the above theorem.

We will attempt to prove that the condition is necessary.
Suppose $\phi\in\Gamma_c(M,E)$ is stationary; that is,
$$
(j^1\phi)^*\inn (j^1Z)\Omega_{\Lag} = 0
\qquad \forall Z\in\vf^{{\rm V}(\pi)}(E)
$$
and consider $X\in\vf (J^1E)$, which can be writen as
$X=X_\phi+X_v$ where $X_\phi$ is tangent to the image of $j^1\phi$ and
$X_v$ is vertical in relation to $\bar\pi^1$,
both in the points of the image of $j^1\phi$.
However, $X_v=(X_v-j^1(\pi^1_*X_v))+j^1(\pi^1_*X_v)$,
where $j^1(\pi^1_*X_v)$ is understood as the prolongation of
a vector field which coincides with $\pi^1_*X_v$ on the image of $\phi$.
Observe that $\pi^1_*(X_v-j^1(\pi^1_*X_v))=0$ on the points
of the image of $j^1\phi$. Therefrom
$$
(j^1\phi)^*\inn (X)\Omega_{\Lag} =
(j^1\phi)^*\inn (X_\phi)\Omega_{\Lag} +
(j^1\phi)^*\inn (X_v-j^1(\pi^1_*X_v))\Omega_{\Lag} +
(j^1\phi)^*\inn (j^1(\pi^1_*X_v))\Omega_{\Lag}
$$
However, $(j^1\phi)^*\inn (X_\phi)\Omega_{\Lag} =0$
because $X_\phi$ is tangent to the image of $j^1\phi$, hence $\Omega_{\Lag}$
acts on linearly dependent vector fields.
Nevertheless, $(j^1\phi)^*\inn (X_v-j^1(\pi^1_*X_v))\Omega_{\Lag} =0$
because $X_v-j^1(\pi^1_*X_v)$ is $\pi^1$-vertical and $\Omega_{\Lag}$
vanishes on these vector fields, when it is restricted to $j^1\phi$
(again, see the comment at the end of the proof).
Therefore
$$
\int_M(j^1\phi)^*\inn (X)\Omega_{\Lag} =
\int_M(j^1\phi)^*\inn (j^1(\pi^1_*X_v))\Omega_{\Lag} =0
$$
since $\phi$ is stationary and $\pi^1_*X_v\in\vf^{{\rm V}(\pi)}(E)$.

\quad\quad ($2\ \Leftrightarrow\ 6$)\quad
$\phi$ is a stationary section if, and only if,
$$
\int_M(j^1\phi)^*\Lie (j^1Z)\Lag =0 \quad ,\quad
\forall Z\in\vf^{{\rm V}(\pi)}(E)
$$
If \dst Z=\beta^A\derpar{}{y^A}\) then
$$
j^1Z=\beta^A\derpar{}{y^A}+\left(\derpar{\beta^A}{x^\mu}+
\derpar{\beta^A}{y^B}v^B_\mu\right)\derpar{}{v_\mu^A}
$$
and hence
$$
\inn (j^1Z)\d\Lag =
\inn (j^1Z)(\d\lag\wedge\omega ) =\inn (j^1Z)(\d\lag\wedge\d^{n+1}\omega ) =
\beta^A\derpar{\lag}{y^A}+\left(\derpar{\beta^A}{x^\mu}+
\derpar{\beta^A}{y^B}v^B_\mu\right)\derpar{\lag}{v_\mu^A}
$$
therefore we have
\beq
(j^1\phi)^*\inn (j^1Z)\d\Lag =
\left[ (\beta^A\circ \phi)\derpar{\lag}{y^A}\circ j^1\phi+
\left(\derpar{\beta^A}{x^\mu}\circ \phi+
\left(\derpar{\beta^A}{y^B}\circ \phi\right)\derpar{\phi^B}{x^\mu}
\derpar{\lag}{v_\mu^A}\circ j^1\phi\right)\right]\d^{n+1}x
\label{ecua1}
\eeq
Now we will calculate the integral. Firstly, one may observe that
the last term gives
$$
\int_M\left(\derpar{\beta^A}{y^B}\circ \phi\right)\derpar{\phi^B}{x^\mu}
\left(\derpar{\lag}{v_\mu^A}\circ j^1\phi\right)\d^{n+1}x=
\int_M\left[\derpar{}{x^\mu}(\beta^A\circ \phi)-
\derpar{\beta^A}{x^\mu}\circ \phi\right]
\left(\derpar{\lag}{v_\mu^A}\circ j^1\phi\right)\d^{n+1}x
$$
and the second term of this last expression cancels the second one in
(\ref{ecua1}).Whereas for the first one, taking into account that
the section is compact supported and integrating by parts, we have
$$
\int_M\derpar{\beta^A\circ \phi}{x^\mu}
\left(\derpar{\lag}{v_\mu^A}\circ j^1\phi\right)\d x^{n+1}x=
-\int_M(\beta^A\circ \phi)
\derpar{}{x^\mu}\left(\derpar{\lag}{v_\mu^A}\circ j^1\phi\right)
\d x^{n+1}x
$$
therefore
$$
\int_M\inn (j^1Z)\d\Lag =\int_M\left(\derpar{\lag}{y^A}\circ j^1\phi-
\derpar{}{x^\mu}\left(\derpar{\lag}{v_\mu^A}\circ j^1\phi\right)\right)
(\beta^A\circ \phi)\d x^{n+1}x=0
$$
and since $\beta^A$ are arbitrary functions, we obtain
$$
\derpar{\lag}{y^A}\circ j^1\phi-
\derpar{}{x^\mu}\left(\derpar{\lag}{v_\mu^A}\circ j^1\phi\right) =0
$$
as set out to prove.
\qed

{\bf Remark}:
\begin{itemize}
\item
In this context, the fundamental theorem of variational calculus is applied
in the following way: let $\phi\in\Gamma_c(M,E)$ be a section for
which there exists $Z\in\vf^{{\rm V}(\pi)}(E)$ and $x\in M$ such that
$$
((j^1\phi)^*\inn (j^1Z)\Omega_{\Lag} )(x)\not= 0
$$
we will see that, in this case, the section is not critical.
In order to do this, it suffices to find a vector field
$Z'\in\vf^{{\rm V}(\pi)}(E)$ such that
$$
\int_M(j^1\phi)^*\inn (j^1Z')\Omega_{\Lag} \not= 0
$$
This vector field $Z'$ is constructed as follows: let $W$ be
a coordinate open set with $x\in W$ and such that
$(j^1\phi)^*\inn (j^1Z)\Omega_{\Lag}\vert_W =g\d^{n+1}x$
with $g>0$, and let $f\colon M\to\Real$ be a function verifying that
\ben
\item
$f\vert_W\leq 0$.
\item
There exists $V\subset W$ such that $f\vert_V>0$.
\item
$f\vert_{M-W}=0$.
\een
Consider the vector field $Z'=fZ$. Since $f$ vanishes out
of the coordinate neighbourhood, the integral involving $Z'$ vanishes there.
Consequently, on $W$, we have
\beann
Z=\beta^A\derpar{}{y^A}
\quad&\Rightarrow&\quad
j^1Z=\beta^A\derpar{}{y^A}+Z=\beta^A\derpar{}{y^A}
\left(\derpar{\beta^A}{x^\mu}+v_\mu^B\derpar{\beta^A}{y^B}\right)
\derpar{}{v_\mu^A}
\\
Z'=fZ
\quad&\Rightarrow&\quad
j^1Z'=j^1(fZ)=fj^1Z+\beta^A\derpar{f}{x^\mu}\derpar{}{v_\mu^A}
\eeann
therefore $\inn (j^1Z')=f\inn (j^1Z)+\inn (X)$,
where \dst X=\beta^A\derpar{f}{x^\mu}\derpar{}{v_\mu^A}\) ,
but since
$$
\Omega_{\Lag} =-\d\derpar{\lag}{v^A_\mu}\wedge\d y^A\wedge\d^n x_\mu
+\d E\wedge\d^{n+1}x
$$
we obtain
\beann
\inn (X)\Omega_{\Lag} &=&
-\beta^B\d\derpar{f}{x^\rho}
\frac{\partial^2\lag}{\partial v^B_\rho\partial v^A_\mu}
\wedge\d y^A\wedge\d^n x_\mu +
\beta^B\derpar{f}{x^\rho}\derpar{E}{v^B_\rho}\wedge\d^{n+1}x
\\
& &-\beta^B\d\derpar{f}{x^\rho}
\frac{\partial^2\lag}{\partial v^B_\rho\partial v^A_\mu}
\wedge\d y^A\wedge\d^n x_\mu +\beta^B\derpar{f}{x^\rho}
\frac{\partial^2\lag}{\partial v^B_\rho\partial v^A_\mu}\wedge\d^{n+1}x
\eeann
If $\phi=(x^\mu,\phi^A(x^\mu ))$, we have
$$
(j^1\phi)\inn (X)\Omega_{\Lag} =-\beta^B\derpar{f}{x^\rho}
\frac{\partial^2\lag}{\partial v^B_\rho\partial v^A_\mu}\Big\vert_{j^1\phi}
\derpar{\phi^A}{x^\mu}\d^{n+1}x +\beta^B\derpar{f}{x^\rho}
\frac{\partial^2\lag}{\partial v^B_\rho\partial v^A_\mu}\Big\vert_{j^1\phi}
\derpar{\phi^A}{x^\mu}\d^{n+1}x=0
$$
therefore
$$
(j^1\phi)\inn (j^1Z')\Omega_{\Lag} =(j^1\phi)f\inn (j^1Z)\Omega_{\Lag}
$$
and, taking into account the conditions on $W$ and $f$, we conclude that
$$
\int_M(j^1\phi)^*\inn (j^1Z')\Omega_{\Lag} =
\int_W(j^1\phi)^*f\inn (j^1Z)\Omega_{\Lag}=
\int_Wf(j^1\phi)^*\inn (j^1Z)\Omega_{\Lag} >0
$$
\end{itemize}

\subsection{Lagrangian formalism with jet fields}

The lagrangian formalism for field theories can also be established in
a more analogous way than in the case of mechanical systems
using {\sl 1-jet fields}. We devote the following paragraps
to developing these topics and review the main results obtained
in the above sections.

First, we will introduce the kind of geometrical objects
which are relevant in order to provide this particular formulation of
dynamics of field theories. Once again, we refer to the appendix
\ref{conec} for the basic concepts related to jet fields
in a jet bundle.

\subsubsection{1-jet fields in the bundle $J^1J^1E$}

Consider the bundle $\bar\pi^1\colon J^1E\to M$.
The jet bundle $J^1J^1E$ is obtained
by defining an equivalence relation on the local sections
of $\bar\pi^1$. Hence, the elements of $J^1J^1E$ are
equivalence classes of these local sections and
$J^1J^1E$ is an affine bundle over $J^1E$,
modelled on the vector bundle
$\bar\pi^{1^*}\Tan^*M \otimes_{J^1E}{\rm V}(\pi^1 )$.
So, we have the commutative diagram
\beq
\begin{array}{ccccc}
J^1J^1E&\mapping{\pi^1_1}&J^1E&\mapping{\pi^1}&E
\\
&
\begin{picture}(20,20)(0,0)
\put(0,0){\mbox{$\bar\pi^1_1$}}
\put(0,20){\vector(1,-1){20}}
\end{picture}
&\Big\downarrow\bar\pi_1&
\begin{picture}(20,20)(0,0)
\put(13,0){\mbox{$\pi$}}
\put(20,20){\vector(-1,-1){20}}
\end{picture}
&
\\
& & M & &
\end{array}
\label{diag}
\eeq

Let ${\cal Y}\colon J^1E\to J^1J^1E$ be a 1-jet field in $J^1J^1E$.
As we know ${\cal Y}$ induces a connection form
$\nabla$ and a horizontal $n+1$-subbundle ${\rm H}(J^1E)$ such that
$$
\Tan J^1E ={\rm V}(\pi^1 )\oplus{\rm Im}\,\nabla =
{\rm V}(\pi^1 )\oplus{\rm H}(J^1E)
$$
where ${\rm H}_{\bar y}={\rm Im}\Tan_{\bar\pi^1 (\bar y)}\psi$,
for $\bar y\in J^1E$
and $\psi\colon M\to J^1E$ a representative of ${\cal Y}(\bar y)$.
We denote by ${\cal D}({\cal Y})$ the $\Cinfty (J^1E)$-module of
sections of ${\rm H}(J^1E)$.

If $(x^\mu ,y^A,v_\mu^A,a_\rho^A,b_{\rho\mu}^A)$ is a natural local
system in $J^1J^1E$, we have the following local expressions
for these elements
\beann
{\cal Y}&=&(x^\mu ,y^A,v_\mu^A,F_\rho^A(x^\mu ,y^A,v_\mu^A),
G_{\nu\rho}^A(x^\mu ,y^A,v_\mu^A))
\\
{\rm H}(J^1E)&=&
\left\{\derpar{}{x^\mu}+F_\mu^A(\bar y)\derpar{}{y^A}+
G_{\rho\mu}^A(\bar y)\derpar{}{v_\rho^A}\right\}
\\
\nabla&=&\d x^\mu\otimes
\left(\derpar{}{x^\mu}+F_\mu^A\derpar{}{y^A}+
G_{\rho\mu}^A\derpar{}{v_\rho^A}\right)
\eeann

If ${\cal Y}\colon J^1E\to J^1J^1E$ is a 1-jet field then
a section $\psi \colon M\to J^1E$ is said to be an
{\sl integral section} of ${\cal Y}$ iff
${\cal Y}\circ \psi =j^1\psi $.
${\cal Y}$ is said to be an {\sl integrable jet field} iff it admits
integral sections.
One may readily check that,
in a natural local system of coordinates in $J^1J^1E$,
$\psi =(x^\mu ,f^\mu (x^\nu ),g^A_\mu (x^\nu )$ and
it is an integral section of ${\cal Y}$ if, and only if,
$\psi $ is a solution of the following system of differential equations
$$
\derpar{f^A}{x^\mu}=F_\mu^A\circ\psi
\qquad \derpar{g_\rho^A}{x^\mu}=G_{\rho\mu}^A\circ\psi
$$
Remember that if $\psi$ is an integral section of ${\cal Y}$,
then the distribution
${\cal D}({\cal Y})$ is tangent to the image of $\psi$ and conversely.
Hence, ${\cal Y}$ is integrable if, and only if,
${\cal D}({\cal Y})$ is an involutive distribution
or, what is equivalent, if, and only if,
the curvature of $\nabla$ is nule.

\subsubsection{The S.O.P.D.E. condition}

The idea of this paragraph is to characterize the integrable jet fields
in $J^1J^1E$ such that their integral sections
are canonical prolongations of sections from $M$ to $E$.

It is well known that there are two natural projections from
$\Tan\Tan Q$ to $\Tan Q$. In the same way, if we
consider the diagram (\ref{diag}), we are going to see that
there is another natural projection from $J^1J^1E$ to $J^1E$.
Let ${\bf y}\in J^1J^1E$ with
\dst{\bf y}\stackrel{\pi^1_1}{\mapsto}\bar y\stackrel{\pi^1}{\mapsto}
y\stackrel{\pi}{\mapsto}x\) ,
and $\psi\colon M\to J^1E$ a representative of ${\bf y}$, that is,
${\bf y} =\Tan_x\psi$. Consider now the section
$\phi =\pi^1\circ\psi\colon M\to E$, then
$j^1\phi (x)\in J^1E$ and we have:

\begin{prop}
The projection
$$
\begin{array}{ccccc}
j^1\pi^1&\colon&J^1J^1E&\longrightarrow&J^1E
\\
& &{\bf y}&\mapsto&j^1(\pi^1\circ\psi )(\bar\pi^1_1({\bf y}))
\end{array}
$$
is differentiable.
\end{prop}
( {\sl Proof} )\quad
Let $(x^\mu ,y^A,v_\mu^A,a_\rho^A,b_{\mu\rho}^A)$ be a natural local
system in a neigbourhood of ${\bf y}_0\in J^1J^1E$ and
${\bf y}_0=(x^\mu_0,y^A_0,{v_\mu^A}_0,{a_\rho^A}_0,{b_{\mu\rho}^A}_0)$.
We have $\pi^1_1({\bf y}_0)=({x^\mu}_0,{y^A}_0,{v_\mu^A}_0)$.
Let $\psi\colon M\to J^1E$ be a representative of ${\bf y}_0$;
locally $\psi=(x^\mu,f^A(x^\nu ),g_\mu^A(x^\nu ))$ with
$$
f^A(x^\nu_0)=y_0^A \quad ,\quad g_\mu^A(x^\nu_0)={v^A_\mu}_0
\quad ;\quad
\derpar{f^A(x^\nu_0)}{x^\rho}={a_\rho^A}_0 \quad ,\quad
\derpar{g_\mu^A(x^\nu_0)}{x^\rho}={b^A_{\rho\mu}}_0
$$
then $\phi =\pi^1\circ\psi =(x^\mu ,f^A(x^\nu ))$, and
$$
j^1\phi =\left( x^\mu ,f^A(x^\nu ),\derpar{f^A}{x^\mu}(x^\nu )\right)
$$
Hence $j^1\pi^1({\bf y}_0)=(x^\mu_0,y^A_0,{a^A_\rho}_0)$
and the result follows.
\qed

{\bf Remark}:
\bit
\item
Observe that $j^1\pi^1$ and $\pi^1_1$ exchange the coordinates
$v^A_\mu$ and $a^A_\mu$.
\eit

\begin{corol}
If $\psi\colon M\to J^1E$ is a section of $\bar\pi^1$,
then $j^1\pi^1\circ j^1\psi =j^1(\pi^1\circ\psi )$.
\end{corol}
\proof
In a local coordinate system $(x^\mu ,y^A,v_\mu^A,a_\rho^A,b_{\mu\rho}^A)$,
we have $\psi =(x^\mu ,f^A(x^\nu ),g_\rho^A(x^\nu ))$ and
\dst j^1\psi =\left( x^\mu ,f^A(x^\nu ),g_\rho^A(x^\nu ),
\derpar{f^A}{x^\rho}(x^\nu ),
\derpar{g_\mu^A}{x^\rho}(x^\nu )\right)\) , but
\dst j^1\pi^1\circ j^1\psi =\left( x^\mu ,f^A(x^\nu),
\derpar{f^A}{x^\rho}(x^\nu )\right) =j^1(\pi^1\circ\psi )\) .
\qed

\begin{definition}
A 1-jet field ${\cal Y}\colon J^1E\to J^1J^1E$
is said to be a
{\rm Second Order Partial Differential Equation (SOPDE)}, or also that
verifies the {\rm SOPDE condition}, iff it
is a section of the projection $j^1\pi^1$ or, what is equivalent,
$$
j^1\pi^1\circ {\cal Y}={\rm Id}_{J^1E}
$$
\end{definition}

This kind of jet fields can be characterized in the following way:

\begin{prop}
Let ${\cal Y}\colon J^1E\to J^1J^1E$ be an integrable 1-jet field.
The necessary and sufficient condition for ${\cal Y}$ to be a SOPDE is that
its integral sections are canonical prolongations of sections
$\phi\colon M\to E$.
\end{prop}
\proof
\quad ($\Longrightarrow$)\quad
If ${\cal Y}$ is a SOPDE then  $j^1\pi^1\circ {\cal Y}={\rm Id}_{J^1E}$.
Let $\psi\colon M\to J^1E$ be an integral section of ${\cal Y}$; that is,
${\cal Y}\circ\psi =j^1\psi$, then
$$
\psi =j^1\pi^1\circ{\cal Y}\circ\psi =
j^1\pi^1\circ j^1\psi =j^1(\pi^1\circ\psi )\equiv j^1\phi
$$
and $\psi$ is a canonical prolongation.

\quad\quad ($\Longleftarrow$)\quad
Now, let ${\cal Y}$ be an integrable jet field whose integral sections
are canonical prolongations. Take $\bar y\in J^1E$ and
$\phi\colon M\to E$ a section such that
$j^1\phi\colon M\to J^1E$ is an integral section of ${\cal Y}$
with $j^1\phi (\pi^1(\bar y))=\bar y$. We have
\beann
(j^1\pi^1\circ{\cal Y})(\bar y)&=&
(j^1\pi^1\circ{\cal Y})(j^1\phi (\bar\pi^1(\bar y)))=
(j^1\pi^1\circ{\cal Y}\circ j^1\phi )(\bar\pi^1(\bar y))=
(j^1\pi^1\circ{\cal Y}\circ\psi )(\bar\pi^1(\bar y))
\\ &=&
(j^1\pi^1\circ j^1\psi)(\bar\pi^1(\bar y))=
j^1(\pi^1\circ\psi)(\bar\pi^1(\bar y))=
j^1\phi (\bar\pi^1(\bar y))=\bar y
\eeann
and ${\cal Y}$ is a SOPDE.
\qed

{\bf Remark}:
\bit
\item
In coordinates, the condition $j^1\pi^1\circ {\cal Y}={\rm Id}_{J^1E}$
is expressed as follows: the jet field
${\cal Y}=(x^\mu ,y^A,v_\mu^A,F_\rho^A,G_{\nu\rho}^A)$
is a SOPDE if, and only if, $F_\rho^A=v_\rho^A$.
\item
If ${\cal Y}=(x^\mu ,y^A,v_\mu^A,v_\rho^A,G_{\nu\rho}^A)$
is a SOPDE and
\dst j^1\phi =\left( x^\mu ,f^A,\derpar{f^A}{x^\mu}\right)\)
is an integral section of ${\cal Y}$, then the differential equations
for $\phi$ are
$$
G_{\nu\rho}^A\left(x^\nu ,f^A,\derpar{f^A}{x^\gamma}\right) =
\frac{\partial^2f^A}{\partial x^\rho\partial x^\nu}
$$
which justifies the nomenclature.
y\eit

\subsubsection{Action of jet fields on forms}

Let ${\cal Y}\colon J^1E\to J^1J^1E$ be a 1-jet field.
A map $\bar{\cal Y}\colon\vf (M)\to\vf (J^1E)$
can be defined in the following way:
Let $Z\in\vf (M)$,
$\bar{\cal Y}(Z)\in\vf(J^1E)$ is the vector field defined as
$$
\bar {\cal Y}(Z)(\bar y):=
(\Tan_{\bar\pi^1 (\bar y)}\psi )(Z_{\bar\pi^1 (\bar y)})
$$
for every $\bar y\in J^1E$ and $\psi\in {\cal Y}(\bar y)$.
This map is an element of $\df^1(M)\otimes_M\vf (J^1E)$
and its local expression is
$$
\bar{\cal Y}\left( f^\mu\derpar{}{x^\mu}\right)=
f^\mu\left(\derpar{}{x^\mu}+F_\mu^A\derpar{}{y^A}+
G_{\rho\mu}^A\derpar{}{v_\rho^A}\right)
$$
(We can recover the connection form using this map,
as follows: $\nabla =\bar{\cal Y}\circ\bar\pi^1_*$).

This map induces an action of ${\cal Y}$ on the forms on
$J^1E$. In fact, let $\xi\in\df^{n+m+1}(J^1E)$,
with $m\geq 0$, we define
$\inn ({\cal Y})\xi\colon\vf (M)\times\stackrel{(n+1)}{\ldots}\times\vf (M)
\longrightarrow\df^m(J^1E)$
given by
$$
((\inn ({\cal Y})\xi )(Z_1,\ldots ,Z_{n+1}))(\bar y;X_1,\ldots ,X_m):=
\xi (\bar y;\bar{\cal Y}(Z_1),\ldots ,\bar{\cal Y}(Z_{n+1}),X_1,\ldots ,X_m)
$$
for $Z_1,\ldots ,Z_{n+1}\in\vf (M)$ and $X_1,\ldots ,X_m\in\vf (J^1E)$.
It is a $\Cinfty (M)$-linear and alternate map on the vector fields
$Z_1,\ldots ,Z_{n+1}$.

\begin{definition}
Let ${\cal Y}\colon J^1E\to J^1J^1E$ be a 1-jet field
and $\xi\in\df^p(J^1E)$.
The $\Cinfty (J^1E)$-linear map $\inn ({\cal Y})\xi$
just defined, extended by zero to forms of degree $p<n+1$,
is called the {\rm inner contraction} of the jet field
${\cal Y}$ and the differential form $\xi$.
\end{definition}

An important result related with the action of integrable jet fields is:

\begin{teor}
Let ${\cal Y}$ be an integrable jet field and $\xi\in\df^{n+2}(J^1E)$.
The necessary and sufficient condition for the integral sections of ${\cal Y}$,
$\psi\colon M\to J^1E$, to verify the relation
$\psi^*\inn (X)\xi =0$, for every $X\in\vf (J^1E)$,
is $\inn ({\cal Y})\xi =0$.
\end{teor}
( {\sl Proof} )\quad ($\Longleftarrow$)\quad
If $\inn ({\cal Y})\xi=0$ then
$((\inn ({\cal Y})\xi )(Z_1,\ldots ,Z_{n+1}))(\bar y;X)=0$,
for every $Z_i\in\vf (M)$, $\bar y\in J^1E$ and $X\in\vf (J^1E)$.
Let $\psi\colon M\to J^1E$ be an integral section of ${\cal Y}$ with
$\psi (x)=\bar y$, then
\beann
0&=&((\inn ({\cal Y})\xi )(Z_1,\ldots ,Z_{n+1}))(\psi (x);X)=
\xi (\psi (x);\Tan_x\psi(Z_1),\ldots ,\Tan_x\psi(Z_{n+1}),X)
\\ &=&
(\inn (X)\xi )(\psi (x);\Tan_x\psi(Z_1),\ldots ,\Tan_x\psi(Z_{n+1}))=
(\psi^*\inn (X)\xi )(x;Z_1,\ldots ,Z_{n+1})
\eeann
for every $X\in\vf (J^1E)$, $x\in M$, $Z_i\in\vf (M)$.
And then the result follows.

\quad\quad ($\Longrightarrow$)\quad
If  $\psi^*\inn (X)\xi =0$, for every $X\in\vf (J^1E)$,
then, let $Z_1,\ldots ,Z_{n+1}\in{\cal X}(M)$, $\bar y\in J^1E$
and $\psi\colon M\to E$ an integral section of ${\cal Y}$
with $\psi (x)=\bar y$. We have
\beann
((\inn ({\cal Y})\xi )(Z_1,\ldots ,Z_{n+1}))(\psi (x);X)&=&
(\inn (X)\xi )(\psi (x);\Tan_x\psi(Z_1),\ldots ,\Tan_x\psi(Z_{n+1}))
\\ &=&
(\psi^*\inn (X)\xi )(x;Z_1,\ldots ,Z_{n+1})=0
\eeann
for every $X\in\vf (J^1E)$, $x\in M$, $Z_i\in\vf (M)$.
Hence $\inn ({\cal Y})\xi =0$.
\qed

 From the definition, one may readily prove the following properties
for this map:

\begin{prop}
\ben
\item
$\inn ({\cal Y})(\xi_1+\xi_2)=\inn ({\cal Y})(\xi_1)+\inn ({\cal Y})(\xi_2)$.
\item
$\inn ({\cal Y})(f\xi )=f\inn ({\cal Y})(\xi )$, for every
$f\in\Cinfty (J^1E)$.
\item
$\inn ({\cal Y})(\xi )\in\df^{n+1}(M)\otimes_M\df^m(J^1E)$.
\een
\end{prop}

The above definition rules also for jet fields in $J^1E$ and differential
forms $\xi\in\df^p(E)$, in an analogous way.
So, every 1-jet field $\Psi$ defines a map
$\bar\Psi\colon\vf (M)\longrightarrow\vf (E)$
as follows: if $y\in E$ then, to every $Z\in\vf (M)$,
$\bar\Psi$ associates the vector field $\bar\Psi (Z)\in\vf (E)$ defined as
$\bar\Psi (Z)(y):=(\Tan_{\pi^1(y)}\phi )(Z_{\pi^1(y)})$
($\phi\in\Psi (y)$).
Then, if $\xi\in\df^{n+m+1}(E)$ ($m\geq 0$), we define
$\inn (\Psi )\xi\colon\vf (M)\times\stackrel{(n+1)}{\ldots}\times\vf (M)
\longrightarrow\df^m(E)$
given by
$$
((\inn (\Psi )\xi )(Z_1,\ldots ,Z_{n+1}))(\bar y;X_1,\ldots ,X_m):=
\xi (\bar y;\bar\Psi (Z_1),\ldots ,\bar\Psi (Z_{n+1}),X_1,\ldots ,X_m)
$$
for $Z_1,\ldots ,Z_{n+1}\in\vf (M)$ and $X_1,\ldots ,X_m\in\vf (E)$.

\subsubsection{Lagrangian and Energy functions. Euler-Lagrange equations}

Now we are ready to establish the lagrangian formulation
of field theories using 1-jet fields in $J^1J^1E$.

First of all, you can check
(using the corresponding local expressions) that the lagrangian
and energy functions can be characterized using jet fields as follows:

\begin{prop}
Let $\ls$ be a lagrangian system and $\omega$ the volume form in $M$. Set
$\moment{Z}{1}{n+1}\in\vf (M)$
such that $\omega (\moment{Z}{1}{n+1})=1$.
\ben
\item
If $\nabla$ is a connection form in $J^1E$,
the lagrangian energy function can be equivalently obtained as:
\ben
\item
For any $1$-jet field $\Psi$ in $J^1E$
$$
\fel =-\Theta_{\Lag}(j^1(\bar\Psi (Z_1),\ldots ,j^1(\bar\Psi (Z_{n+1}))
$$
\item
For any 1-jet vector ${\cal Y}$ in $J^1J^1E$,
$$
\fel =(\inn ({\cal Y})\del )(\moment{Z}{1}{n+1})
$$
\een
\item
For any 1-jet field ${\cal Y}$ in $J^1J^1E$,
the lagrangian function can be recover as
$$
\lag =(\inn ({\cal Y})\Lag )(\moment{Z}{1}{n+1})
$$
\een
\end{prop}

On the other hand, it is known that the critical section
of the associated variational problem satisfies the following conditions:
\ben
\item
They are sections $\phi\colon M\to E$.
\item
$(j^1\phi )^*\inn (j^1Z)\Omega_{\Lag}=0$, for every $Z\in\vf (J^1E)$.
\een
then, taking into account the statements and comments
in the above sections, it is evident that:

\begin{prop}
Let $\ls$ be a lagrangian system.
The critical section of the variational problem posed by $\Lag$
are the integral sections of a jet field ${\cal Y}$, if and only if, ${\cal Y}$
satisfies the following conditions:
\ben
\item
${\cal Y}$ is integrable; that is, ${\cal R}({\cal Y})=0$.
\item
${\cal Y}$ is a SOPDE; that is, $j^1\pi^1\circ {\cal Y}={\rm Id}_{J^1E}$.
\item
$\inn ({\cal Y})\Omega_{\Lag}=0$.
\een
This is the version of the Euler-Lagrange equations in terms of jet fields.
\end{prop}

\section{Lagrangian symmetries and Noether's theorem}

\subsection{Symmetries. Noether's theorem}

We are going to consider the question of symmetries in the
lagrangian formalism of field theories.
The most common and meaningful way to treat this subject is
by using the theory of {\sl actions of Lie groups and algebras}.
However, in this chapter, we will give a brief and simple
introduction to this question and we reserve the study
of symmetries as actions of Lie groups for further research.

\subsubsection{Symmetries of lagrangian systems}

In Mechanics, a {\sl symmetry} of a lagrangian dynamical system
is a diffeomorphism in the phase space of the system
(the {\sl tangent bundle}) which leaves the lagrangian function invariant.
Thus, it is usual to speak about {\sl natural} or {\sl physical symmetries},
that is, those which are canonical liftings of diffeomorphisms in the
basis of the tangent bundle: the {\sl configuration space} of the system
(or also vertical diffeomorphisms).
These preserve the canonical structures of the tangent bundle and, also,
as a consequence, the dynamical geometric structures induced
by the lagrangian.

Nevertheless, diffeomorphisms in the tangent bundle
which leave the lagrangian invariant, but which are not canonical liftings
of diffeomorphisms in the basis, do not preserve either
the natural or the dynamical structures of the tangent bundle.
Thus, in the lagrangian formalism,
these kinds of symmetries are usually considered as pathological
or undesired from the physical point of view and in general
they are not considered (although the situation is not
the same in the hamiltonian formalism, where symmetries which are not
canonical liftings play a relevant role in some theories).

Taking these comments as our standpoint, and in order to
establish a more general concept of symmetry for field theories, we define:

\begin{definition}
Let $\ls$ be a lagrangian system.
\ben
\item
A {\rm natural symmetry} of the lagrangian system is a diffeomorphism
$\Phi\colon E\to E$ such that its canonical prolongation
$j^1\Phi\colon J^1E\to J^1E$ leaves $\Lag$ invariant.
\item
A {\rm (general) symmetry} of the lagrangian system is a diffeomorphism
$\Upsilon\colon J^1E\to J^1E$ such that
\ben
\item
$\Upsilon$ leaves the canonical geometric structures of $J^1E$
($\theta$, ${\cal M}_c$, {\cal S}, ${\cal V}$) invariant.
\item
$\Upsilon$ leaves $\Lag$ invariant.
\een
\een
\label{sym}
\end{definition}

It is evident that every natural symmetry is also a (general) symmetry,
but the converse is not true.
At present, we are interested only in natural symmetries,
although all the comments and results we obtain
for them also hold for (general) symmetries.

{\bf Remarks}:
\begin{itemize}
\item
Since $\Lag$ is a $(n+1)$-form on $J^1E$, $j^1\Phi$ acts on it by pull-back.
Then, for $\Lag$ to be invariant by $j^1\Phi$ means that
$(j^1\Phi )^*\Lag =\Lag$. Hence,
one way to obtain the explicit expression of $(j^1\Phi )^*\Lag$
consists in considering the diagram
$$
\begin{array}{ccccc}
& (\vf (J^1E))^{n+1} & $\rightarrowfill$ &
(\vf (J^1E))^{n+1} &
\\
& & ((j^1\Phi )_*)^{n+1} & &
\\
(j^1\Phi )^*\Lag & \Big\downarrow & & \Big\downarrow & \Lag
\\
& & ((j^1\Phi )^*)^{-1} & &
\\
& \Cinfty (J^1E) & $\rightarrowfill$ & \Cinfty (J^1E) &
\end{array}
$$
where $\Lag$ is understood as an element of $\df^{n+1}(J^1E)$,
we have
$$
(j^1\Phi )^*\Lag = (j^1\Phi )^*\circ\Lag\circ ((j^1\Phi )_*)^{n+1}
$$
that is, if $X_i\in\vf (J^1E)$ and $\bar y\in J^1E$,
$$
((j^1\Phi )^*\Lag )(\bar y;\moment{X}{1}{{n+1}})=
\Lag (j^1\Phi (\bar y);(j^1\Phi )_*X_1,\ldots ,(j^1\Phi )_*X_{n+1})
$$
\item
Another way to obtain explicit expression of $(j^1\Phi)^*\Lag$
consists in considering $\Lag$ as an element of
$\Gamma (J^1E,\bar\pi^{1^*}\Lambda^{n+1}\Tan^*M,)$,
then we have the diagram
$$
\begin{array}{ccccc}
& \bar\pi^{1^*}\Lambda^{n+1}\Tan^*M &
$\rightarrowfill$ &
\bar\pi^{1^*}\Lambda^{n+1}\Tan^*M &
\\
& & (j^1\Phi ,\Lambda^{n+1}{\Phi_M^*}^{-1}) & &
\\
(j^1\Phi )^*\Lag & \Big\downarrow & & \Big\downarrow & \Lag
\\
& & j^1\Phi & &
\\
& J^1E & $\rightarrowfill$ & J^1E &
\end{array}
$$
and then we obtain that
$$
(j^1\Phi )^*\Lag =
(j^1\Phi ,\Lambda^{n+1}{\Phi_M^*}^{-1})^{-1}\circ\Lag\circ j^1\Phi =
(j^1\Phi^{-1} ,\Lambda^{n+1}\Phi_M^*)\circ\Lag\circ j^1\Phi
$$
that is, the same expression as above, taking into account that $\Lag$ is
a $(n+1)$-form on $M$ with coefficients in $\Cinfty (J^1E)$.
\item
Finally, for the last interpretation we consider the diagram
$$
\begin{array}{ccccc}
& (\vf (E))^{n+1} & $\rightarrowfill$ & (\vf (E))^{n+1} &
\\
& & (\Phi_M)_*^{n+1} & &
\\
(j^1\Phi )^*\Lag & \Big\downarrow & & \Big\downarrow & \Lag
\\
& & ((j^1\Phi )^*)^{-1} & &
\\
& \Cinfty (J^1E) & $\rightarrowfill$ & \Cinfty (J^1E) &
\end{array}
$$
and therefore
$$
(j^1\Phi )^*\Lag = (j^1\Phi )\circ\Lag\circ (\Phi_M)_*^{n+1}
$$
\end{itemize}

As predicted, we have the following relevant property:

\begin{prop}
If $\Phi\colon E\to E$ is a natural symmetry of the lagrangian system
$\ls$ then
\ben
\item
$(j^1\Phi )^*\vartheta_{\Lag} =\vartheta_{\Lag}$.
\item
$(j^1\Phi )^*\Theta_{\Lag} =\Theta_{\Lag}$.
\een
\label{invasim}
\end{prop}
( {\sl Proof} )\quad
Since $j^1\Phi\colon J^1E\to J^1E$ is a natural prolongation, it
verifies that
$$
(j^1\Phi )^*\theta =\theta \qquad , \qquad
(j^1\Phi )^*{\cal S} ={\cal S} \qquad , \qquad
(j^1\Phi )^*\Theta_{\Lag} =\Theta_{\Lag}
$$
and, since the pull-back commutes with the contractions and
$(j^1\Phi )^*\d\Lag =\d\Lag$, the desired results follow immediately.
\qed

\subsubsection{Infinitesimal symmetries}

A symmetry of a lagrangian system can be thought of being
generated by a vector fied.
This is the standpoint for the following definition:

\begin{definition}
Let $\ls$ be a lagrangian system.
\ben
\item
An {\rm infinitesimal natural symmetry} of the lagrangian system is a
vector field $Z\in\vf (E)$ such that its canonical prolongation
leaves $\Lag$ invariant:
$$
\Lie (j^1Z)\Lag =0
$$
\item
An {\rm infinitesimal (general) symmetry} of the lagrangian system is a
vector field $X\in\vf (J^1E)$ such that:
\ben
\item
The canonical geometric structures of $J^1E$
($\theta$, ${\cal M}_c$, {\cal S}, ${\cal V}$) are invariant
under the action of $X$.
\item
$X$ leaves $\Lag$ invariant.
\een
\een
In these cases $\Lag$ is said to be {\rm invariant} by $j^1Z$ or
$X$ respectively.
\label{insym}
\end{definition}

As in the above paragraph we will limit ourselves to natural symmetries.

{\bf Remarks}:
\begin{itemize}
\item
$\Lie (j^1Z)\Lag$ is calculated in the usual way, once $\Lag$ is understood
as a $(n+1)$-form on $J^1E$; that is
$$
\Lie (j^1Z)\Lag =\inn (j^1Z)\d\Lag +\d\inn (j^1Z)\Lag
$$
\item
The statement that $Z\in\vf (E)$ is an infinitesimal symmetry is equivalent
to stating that the local one-parameter groups of $Z$ are symmetries of the
lagrangian system (in the sense of definition \ref{sym}).
\item
This kind of symmetries, as well as those in the above paragraph,
are said to be {\sl dynamical symmetries} (in order to distinguish them
from other kinds of symmetries), because they are canonical
prolongations of diffeomorphisms or vector fields on $E$, which is the
configuration space of the physical fields.
\end{itemize}

In order to obtain the local coordinate expression, suppose that
\beann
Z&=&\alpha^\mu\derpar{}{x^\mu}+\beta^A\derpar{}{y^A}
\\
j^1Z&=&\alpha^\mu\derpar{}{x^\mu}+\beta^A\derpar{}{y^A}+
\left(\derpar{\beta^A}{x^\mu}-\derpar{\alpha^\rho}{x^\mu}v^A_\rho+
\derpar{\beta^A}{y^B}v^B_\mu\right)\derpar{}{v_\mu^A}
\eeann
and, taking $\Lag =\lag\omega$, we have
$\d\Lag =\d\lag\wedge\omega$, hence
\beann
\inn (j^1Z)\d\Lag &=&((j^1Z)\lag )\omega -\d\lag\wedge\inn (j^1Z)\omega
\\
\inn (j^1Z)\Lag &=&\lag\inn (j^1Z)\omega
\\
\d\inn (j^1Z)\Lag &=&\d\lag\wedge\inn (j^1Z)\omega +\lag\d\inn (j^1Z)\omega
\eeann
therefore
\beann
\Lie (j^1Z)\Lag &=&((j^1Z)\lag )\omega +\lag\d\inn (j^1Z)\omega=
((j^1Z)\lag )\omega
+\lag\d\left(\alpha^\mu\inn\left(\derpar{}{x^\mu}\right)\omega\right)
\\ &=&
((j^1Z)\lag )\omega +\lag\sum_\mu\derpar{\alpha^\mu}{x^\mu}\omega=
\left((j^1Z)\lag +\lag\sum_\mu\derpar{\alpha^\mu}{x^\mu}\right)\omega
\eeann
therefore, if $\Lag$ is invariant by $Z$ we obtain
$$
\delta_Z\lag\equiv (j^1Z)\lag +\lag\sum_\mu\derpar{\alpha^\mu}{x^\mu} = 0
$$
where $\delta_Z\lag$ is usually called the {\sl total variation} of $\lag$.

In the particular case that $Z$ is a vertical vector field,
$\alpha^\mu =0$ and $\Lag$ is invariant by $Z$ if, and only if,
$(j^1Z)\lag =0$.

In relation to infinitesimal symmetries,
we have the following relevant property:

\begin{prop}
If $Z\in\vf (E)$ is an infinitesimal natural symmetry
of the lagrangian system $\ls$ then:
\ben
\item
$\Lie (j^1Z)^*\vartheta_{\Lag} =0$.
\item
$\Lie (j^1Z)^*\Theta_{\Lag} =0$.
\een
\end{prop}
( {\sl Proof} )\quad
It is a direct consequence of proposition \ref{invasim}.
\qed

\subsubsection{Noether's theorem}

Closely related  to the problem of symmetries is the question of
{\sl conserved quantities}. The main result is:

\begin{teor}{\rm (Noether):}
Let $Z$ be an infinitesimal symmetry of the  system
$\ls$. Then, the $n$-form ${\rm J}(Z):=\inn (j^1Z)\Theta_{\Lag}$
is constant on the critical sections of the variational problem
posed by $\Lag$.
\end{teor}
( {\sl Proof} )\quad
Let $\phi\colon M\to E$ be a critical section of the variational problem,
that is,
$$
(j^1\phi)^*\inn (X)\d\Theta_{\Lag} =0 \quad ,\quad \forall X\in\vf (J^1E)
$$
Since $\Lag$ is invariant under $Z$ we have
$$
0=\Lie (j^1Z)\Theta_{\Lag} =
\d\inn (j^1Z)\Theta_{\Lag} +\inn (j^1Z)\d\Theta_{\Lag}
$$
therefore
\beann
0=(j^1\phi)^*[\Lie (j^1Z)\Theta_{\Lag} ]&=&
(j^1\phi)^*[\d\inn (j^1Z)\Theta_{\Lag} ]
+(j^1\phi)^*[\inn (j^1Z)\d\Theta_{\Lag} ]
\\ &=&
(j^1\phi)^*[\d\inn (j^1Z)\Theta_{\Lag} ]=
\d [(j^1\phi)^*\inn (j^1Z)\Theta_{\Lag} ]
\eeann
and the resulst follows.
\qed

\begin{definition}
Let $\ls$ be a lagrangian system.
If $\phi$ is a critical section for the variational problem posed by $\Lag$,
then the expression $(j^1\phi)^*{\rm J}(Z)$ is called
{\sl Noether's current} associated with $Z$.
\end{definition}

In order to give an interpretation of this current,
let $N\subset M$ be a manifold with boundary $\partial N$,
such that ${\rm dim}\, N={\rm dim}\, M$ and
$\phi$ a critical section. Since
$\d [(j^1\phi)^*\inn (j^1Z)\Theta_{\Lag} ]=0$, we obtain
$$
0=\int_N\d [(j^1\phi)^*\inn (j^1Z)\Theta_{\Lag} ]=
\int_N\d [(j^1\phi)^*\inn (j^1Z)\Theta_{\Lag} ]=
\int_{\partial N}(j^1\phi)^*\inn (j^1Z)\Theta_{\Lag}
$$

\subsubsection{Noether's theorem for jet fields}

Next we are going to interpret the idea of {\sl conserved quantity}
in terms of jet fields.
In analogy with time-dependent lagrangian mechanical systems
and taking into account how a jet field ${\cal Y}$ acts on the differential
forms, it seems reasonable to find invariants which are
differential $n$-forms $\zeta\in\df^n(J^1E)$ such that
$\inn ({\cal Y})\d\zeta =0$.

Then, in order to state Noether's theorem in this context,
we must first recall the concepts of natural symmetry and
infinitesimal natural symmetry
given in definitions \ref{sym} and \ref{insym} and their properties.
At this point, the following result is required:

\begin{lem}
Let ${\cal Y}\colon J^1E\to J^1J^1E$ be an integrable jet field such that
it is a SOPDE. If $\zeta\in\df^{n+1}(J^1E)$ belongs to the ideal
generated by the contact module,
(which is denoted by ${\cal I}({\cal M}_c)$), then
$$
\inn ({\cal Y})\zeta =0
$$
\end{lem}
( {\sl Proof} )\quad
Consider $\bar y\in J^1E$ and let $\phi\colon M\to E$ be a section
such that $j^1\phi$ is an integral section of ${\cal Y}$ passing
through $\bar y$. In this case we have the map
$$
\inn ({\cal Y})\zeta\colon \vf (M)\times\stackrel{n+1}{\ldots}\vf (M)
\longrightarrow \Cinfty (J^1E)
$$
such that
$$
(\inn ({\cal Y} )\zeta )(\bar y;Z_1,\ldots ,Z_{n+1})=
\zeta (\bar y;(j^1\phi )_*Z_1,\ldots ,(j^1\phi )_*Z_{n+1})=
((j^1\phi )^*\zeta )(\bar\pi^1 (\bar y);Z_1,\ldots ,Z_{n+1}) = 0
$$
because $\zeta$ belongs to the ideal generated by the contact module and
then $(j^1\phi )^*\zeta =0$.
\qed

Now we are ready to prove that:

\begin{teor} {\rm (Noether)}:
Let $\ls$ be a lagrangian system and
${\cal Y}\colon J^1E\to J^1J^1E$ a jet field verifying the
Euler-Lagrange equations for the lagrangian density $\Lag$.
Let $X\in\vf (J^1E)$ be a vector field satisfying the following conditions:
\ben
\item
$X$ preserves the contact module:
$\Lie (X){\cal M}_c\subset{\cal M}_c$.
\item
There exist $\xi\in\df^n(J^1E)$, and
$\alpha\in{\cal I}({\cal M}_c)\subset\df^{n+1}(J^1E)$,
such that
$\Lie (X)\Lag =\d\xi +\alpha$.
\een
Then the following relation holds
$$
\inn ({\cal Y})\d (\xi -\inn (X)\Theta_{\Lag}) =0
$$
\end{teor}
( {\sl Proof} )\quad
On the one hand we have
$$
\Lie (X)\Theta_{\Lag} =\inn (X)\d\Theta_{\Lag} +\d\inn (X)\Theta_{\Lag}
$$
however
$$
\Lie (X)\Theta_{\Lag} =\Lie (X)(\vartheta_{\Lag} +\Lag ) :=
\beta +\d\xi +\alpha
$$
where $\beta\in{\cal I}({\cal M}_c)$ because
$\vartheta_{\Lag}\in{\cal I}({\cal M}_c)$. Therefrom
$$
\d (\xi -\inn (X)\Theta_{\Lag})=\inn (X)\d\Theta_{\Lag}+\zeta
$$
with $\zeta\in{\cal I}({\cal M}_c)$. Hence
$$
\inn ({\cal Y})\d (\xi -\inn (X)\Theta_{\Lag})=0
$$
since $\inn ({\cal Y} )\inn (X)\d\Theta_{\Lag}=0$,
for every $X\in\vf (J^1E)$,  because ${\cal Y}$ is a solution of the
Euler-Lagrange equations, and
$\inn ({\cal Y})\zeta =0$ as a consequence of the last lemma.
\qed

\section{Examples}

\subsection{Electromagnetic field (with fixed background)}

In this case $M$ is space-time endowed with a semi-riemannian metric $g$,
$E=\Tan^*M$ is a vector bundle over $M$ and
$\pi_M\colon \Tan^*M\to M$ denotes the natural projection.
Sections of $\pi_M$ are the so-called {\sl electromagnetic potentials}.
Using the linear connection associated with the metric $g$,
one can be sure that $J^1E\to \Tan^*M$
is a vector bundle. It then coincides with its associated
vector bundle, and, since ${\rm V}E=\pi_M ^*\Tan^*M$, we have
$J^1E=\pi_M ^*\Tan^*M\otimes_E \pi_M ^*\Tan^*M$.

Coordinates in $J^1E$ are usually denoted $(x^\mu, A_\nu,v_{\mu\nu})$
\footnote{
Sometimes, the notation $A_{\nu ,\mu}$ is used instead of $v_{\mu\nu}$.
}.
The coordinates $A_1,A_2,A_3$ constitute the {\sl vector potential}
and $A_4$ is the {\sl scalar potential},
if $M=\Real^4$ and the metric is $+++-$.

Let $\phi\colon M\to \Tan^*M$ be a section of $\pi$.
Then $j^1\phi\colon M\to\pi^*\Tan^*M\times\pi^*\Tan^*M$ is just
$j^1\phi =\Tan\phi$. Locally $\phi =\phi_\nu\d x^\nu$, therefore
\dst j^1\phi =\derpar{\phi_\nu}{x^\mu}\d x^\mu\otimes\d x^\nu\) .
Then, one may observe that $j^1\phi$ is a metric tensor on $M$.

The lagrangian density is the following:
consider $\bar y\in J^1E$ and let $\phi\colon M\to\Tan^*M$
be a representative of $\bar y$; that is, $\bar y=\Tan_{\bar\pi (\bar y)}\phi$,
then
$$
\Lag =\frac{1}{2}\|\d\phi\|^2\d V_g
$$
where $\d\phi$ is the exterior differential of $\phi$,
$\|\ \|$ is the norm induced by the metric $g$ on the
$2$-forms on $M$ and $\d V_g$ is the volume element
associated with the metric $g$.
Observe that $\d\phi$ is the skew-symmetric part of the matrix
giving $\Tan\phi$ or, in other words, the skew-symmetric part of
the metric $\Tan\phi$ on $M$.
It is usual to write $\phi$ in the form
$A=A_\mu\d x^\mu$ and $\d A=A_{\mu ,\nu}\d x^\nu\wedge\d x^\mu$,
as we stated above.

The symmetries of the problem arise from the isometries of $g$ and from
the fact that the symmetric part of $\Tan\phi$ has not been taken into account
in order to construct the lagrangian density.

Finally, if $f\in\Cinfty (M)$, then $f$ acts on $\Tan^*M$
in the following way:
$$
\begin{array}{ccccc}
\tilde f & \colon &\Tan^*M & \longrightarrow & \Tan^*M
\\
& & (x,\alpha ) & \mapsto & (x,\alpha +\d f (x))
\end{array}
$$
$\tilde f$ being a diffeomorphism which induces the identity on $M$.
Moreover, if $\tilde y\in J^1E$, the equality
$\Lag (j^1\tilde f(\tilde y))=\Lag (\tilde y)$
holds trivially.
This is the so-called {\sl gauge invariance}
of the electromagnetic field.

\subsection{Bosonic string}

Let $M$ be a two-dimensional manifold and $B$ a $d+1$-dimensional manifold
endowed with a given metric $g$, with signature $(+,\ldots,+,-)$,
that is, a  Lorentz metric.
A {\sl bosonic string} in $B$ is a map $\phi\colon M\longrightarrow B$.

Furthermore, in order to construct the lagrangian function we need a
metrics on $M$, we take the bundle
$$
E=M\times B\times{\rm S}^{1,1}_2(M)
$$
where ${\rm S}^{1,1}_2(M)$ is the  bundle over $M$ whose
sections are the $2$-covariant tensors with signature $(+,-)$.
In this case we have:
$$
J^1E = J^1(M\times B)\times J^1({\rm S}^{1,1}_2(M))
$$
and we have obtained the so-called {\it Polyakov approach}
where the metrics on $M$ is dynamics.

In this approach, a field $\psi$ is a couple $(\phi ,h)$,
where $\phi\colon M\to B$ is a differentiable map (the {\sl string})
and $h$ is a Lorentz metric in $M$.
If $\bar y\in J^1E$ and $\psi =(\phi ,h)$ is a
representative of $\bar y$; that is, $(j^1\psi )(\bar\pi^1(\bar y))=\bar y$,
the lagrangian density at the point $\bar y$ is given by
$$
\Lag (\bar y)=-\frac{1}{2}\inn (h^*)\phi^*g\d V_h
$$
where $h^*$ is the dual metric of $h$,
$\inn (h^*)\phi^*g$ denotes the total contraction of the
$2$-covariant tensor $\phi^*g$ with the $2$-contravariant tensor $h^*$
and $\d V_h$ is the volume element associated with the metric $h$.

Denoting by $x^\mu$ the coordinates on $M$, by $y^A$ the ones in $B$
and by $h_{\sigma\rho}$ the coordinates in the fiber of
${\rm S}^{1,1}_2(M)$, then the coordinates in $J^1E$ are
$(x^\mu, y^A,h_{\sigma\rho},v^A_\mu,w_{\sigma \rho \mu})$,
$\mu=1,2,\ A=1,\ldots,d,0$.
In such a local system, if $\phi =\phi^A$, $g=g_{AB}\d y^A\otimes\d y^B$
$h=h_{\sigma\rho}\d x^\sigma\otimes\d x^\rho$
and \dst h^*=h^{\sigma\rho}\derpar{}{x^\sigma}\otimes\derpar{}{x^\rho}\) ,
the linear expression of $\Lag$ is
$$
\Lag =
-\frac{1}{2}h^{\sigma\rho}g_{AB}v_\sigma^Av_\rho^B\sqrt{|{\rm det}\, h|}\d^2x
$$

In the {\it Nambu approach} we use the metric
induced on $M$ by the string itself; that is, $h=\phi^*g$,
where $g$ is the metric over $B$. In this  case,
$E=M\times B$ and the rest are as in the Polyakov approach
with the appropriate changes in the coordinates.

In this last case, a field is a map $\phi\colon M\to B$.
Then, if $\bar y\in J^1E$ and $\phi$ is a representative of it,
the lagrangian density is
$$
\Lag (\bar y)=\frac{1}{2}\d V_{\phi^*g}
$$
A local system of canonical coordinates is now
$(x^\mu, y^A,v^A_\mu)$ and if we write $\phi =\phi^A$
and $g=g_{AB}\d y^A\otimes\d y^B$, we have
$$
\Lag =\frac{1}{2}
\sqrt{\Big |{\rm det}\, \left(g_{AB}\derpar{\phi^A}{x^1}\derpar{\phi^B}{x^2}
\right)\Big |}\d x^1\d x^2
$$
thus it is easy to see that this expression is the infinitesimal ``area''
of the parallelogram generated by
\dst\Tan\phi\left(\derpar{}{x^1}\right)\)
and \dst\Tan\phi\left(\derpar{}{x^2}\right)\) in $B$.

\subsection{WZWN model}

Consider $M=S^2$ endowed with a riemannian metric $h$.
Let $G$ be a Lie group, then take $E=M\times G$
as the manifold of configuration. If ${\bf g}$ is the Lie algebra of $G$,
let $\Omega \in\df^1(G)\otimes{\bf g}$ be the Maurer-Cartan form
of $G$ and denote $K$ the Killing metrics on ${\bf g}$.
A section of $E$ is given by a map $g\colon M\to G$ and defines the function
$(h^*\otimes K)(g^*\Omega, g^*\Omega )\in\Cinfty (M)$,
$h^*$ being the dual metrics of $h$.
This function is used in order to define another one
in the $1$-jet bundle $J^1E=J^1(M\times G)$ in a natural way
as follows: if $\bar y\in J^1E$, with $\bar\pi^1 (\bar y)=x$,
and $g\colon M\to G$ is a representative of it,
we have that the map
$$
x \longmapsto (h^*\otimes K)(g_x^*\Omega, g_x^*\Omega )
$$
does not depend on the chosen representative and hence we can define a function
$$
\bar y \longmapsto
(h^*\otimes K)(g^*_{\bar\pi^1 (\bar y)}\Omega ,g^*_{\bar\pi^1 (\bar y)}\Omega )
$$

The first approach to the WZWN theory takes ${\bf SU}(2)$
as the Lie group, that is
$$
{\bf SU}(2)={\bf SL}(2,\Complex )\cap{\bf U}(2)=
\{ A\in{\bf GL}(2,\Complex )\ |\ AA^\dagger =I,\ {\rm det}\, A=1\}
$$
(Observe that ${\bf SU}(2)$ is diffeomorphic to $S^3$).
In this case the fields are maps $g\colon S^2\to {\bf SU}(2)$.

The second approach consists in considering that the metrics $h$ on
$S^2$ is dynamical too. In that case the configuration bundle is
$E'=(S^2\times G)\times_{S^2}{\cal S}_2(S^2)$,
where ${\cal S}_2(S^2)$ is the bundle over $S^2$ whose sections are
the riemannian metrics on $S^2$.
A section of $E'$ is given by a pair $(g,h)$, where
$g\colon S^2\to G$ is a map and $h$ is a riemannian metrics on $S^2$.
The $1$-jet bundle is $J^1E'=J^1(S^2\times G)\times J^1({\cal S}_2(S^2))$.

Like in the first case we have a function defined
in a natural way in $J^1E'$. In fact, let $\bar y\in J^1E'$
and let $(g,h)$ be a representative of $\bar y$. The map
$$
\bar y \longmapsto (h^*_{\bar\pi^1 (\bar y)}\otimes K)
(g^*_{\bar\pi^1 (\bar y)}\Omega ,g^*_{\bar\pi^1 (\bar y)}\Omega)
$$
is independent of the representative. We use this function
in order to construct the lagrangian density of this model which is
$$
\Lag (\bar y) :=
\bar y \longmapsto (h^*_{\bar\pi^1 (\bar y)}\otimes K)
(g^*_{\bar\pi^1 (\bar y)}\Omega ,g^*_{\bar\pi^1 (\bar y)}\Omega) \d V_h
$$
(This corresponds only to the so-called
{\sl non topological} or {\sl dynamical term} of the lagrangian of the
WZWN model).

The systems of natural coordinates are
$(x^\mu,g^A,h_{\rho\sigma},v^A_\mu ,h_{\rho\sigma\mu})$,
where $x^\mu$ are coordinates in $S^2$ ($\mu =1,2$),
$g^A$ are coordinates in $G$ and $h_{\rho\sigma}$
makes up a local basis of sections of ${\cal S}_2(S^2)$.

{\bf A comment about ${\bf SU}(2)$}:
\begin{itemize}
\item
${\bf SU}(2)$ is the universwal covering of ${\bf SO}(3)$,
the group of rotations in the euclidean space of dimension $3$.
Thus, ${\bf SU}(2)$ is simply connected. An element of this group
is a matrix $A$ as the following
$$
A=\left(\matrix{\alpha & \beta \cr -\bar\beta & \bar\alpha \cr}\right)
\quad\quad
\|\alpha\|^2+\|\beta\|^2=1
\quad\quad
\alpha ,\beta\in\Complex
$$
Therefrom, it is obvious that ${\bf SU}(2)$ is diffeomorphic to
$S^3\subset\Real^4$.

The Lie algebra ${\bf su}(2)$ of ${\bf SU}(2)$ is generated by
$$
E_1=\frac{1}{2}\left(\matrix{ 0 & -i \cr -i & 0 \cr}\right)
\quad ; \quad
E_2=\frac{1}{2}\left(\matrix{ 0 & -1 \cr 1 & 0 \cr}\right)
\quad ; \quad
E_3=\frac{1}{2}\left(\matrix{ -i & 0 \cr 0 & i \cr}\right)
$$
Observe that $E_j=-\frac{1}{2}i\sigma_j$, $\sigma_j$ being
the Pauli matrices. If these matrices $\sigma_j$ are used as a basis of
${\bf SU}(2)$ then the coefficients are imaginary numbers.
\end{itemize}

\section{Conclusions and outlook}

In this work we have formulated a purely geometric theory of the
lagrangian formalism of first-order classical field theories.
The main goals we have achieved are the following:
\ben
\item
We give an intrinsic definition of all the geometric elements
associated with $1$-jet bundles:
the {\sl vertical diferential}, the {\sl canonical form},
the {\sl contact module} and the {\sl vertical endomorphisms}.
Some of them,
namely the structure canonical form, have been defined previously
\cite{Gc-74}; but in other cases, such as the vertical endomorphisms,
the construction presented here is original.
In addition, we have given a detailed description of
the  {\sl canonical prolongations} of sections,
difeomorphisms and vector fields, as well as the properties
of invariance of all the above defined geometric elements
by these prolongations.
\item
Using some of the previously introduced geometric elements,
we have constructed the {\sl lagrangian forms} in an intrinsic way,
with no use of connections \cite{Gc-74}.
Neither is the Legendre transformation used
to obtain lagrangian forms from the dual (hamiltonian) formalism of the theory,
as is the case in other works \cite{GS-73}, \cite{GIMMSY-mm}.
\item
Nevertheless, this is not the case when we define
the {\sl density of lagrangian energy} and the {\sl lagrangian energy function}
in an intrinsic way since, as we have specified,
such a construction necessarily requires the use of a connection.
We also need this connection when we wish to dualize the theory;
that is, in order to obtain the {\sl hamiltonian formalism},
as will be seen in a forthcoming work \cite{EMR-ghft}.
\item
After analyzing the {\sl variational problem} obtaining the
{\sl critical sections}, another outstanding point is that
we achieve a unified presentation of the field equations
in field theories and the dynamical equations for
non-autonomous mechanical systems, by means of the
{\sl $1$-jet field formulation}. The problem of the integrability
of the evolution equations and the fact that these equations are
second-order partial differential equations is clearly interpreted
in a geometrical way using this formalism.
\item
Furthermore, we have carried out an introduction to the study of
{\sl lagrangian symmetries}, classifying them in several ways and
stating and proving a geometrical version of
{\sl Noether's theorem} based on the $1$-jet field formalism.
\item
As a final remark, we have briefly dealt with some classical examples,
making the geometrical structure of each one of them evident.
\een

Among the problems we will study in forthcoming papers,
in which we hope to provide an original treatment,
we draw attention to the following:
\bit
\item
How to construct the hamiltonian formalism of
first-order classical field theories.
\item
An investigation of {\sl non-regular} classical field theories.
\item
To make a detailed analysis of symmetries in field theories.
\eit

\appendix

\section{Connections and jet fields in a first-order jet bundle}
\protect\label{conec}

In this appendix, we summarize some elements of the theory of
jet fields and connections
in fiber bundles (see \cite{Sa-89} for details).

\subsection{Basic definition and properties}

In order to set the main definition of this appendix,
first we prove the following statement:

\begin{prop}
Let $\pi\colon E\to M$ be a fiber bundle and
$\pi^1\colon J^1E\to E$ the corresponding first-order jet bundle.
The following elements can be canonically constructed one from the other:
\ben
\item
A $\pi$-semibasic $1$-form $\nabla$ on $E$
with values in $\Tan E$ (that is, an element of
$\Gamma (E,\pi^*\Tan^*M)\otimes\Gamma (E,\Tan E)$),
such that $\nabla^*\alpha =\alpha$, for every
$\pi^1$-semibasic form $\alpha\in\df^1(E)$.
\item
A subbundle ${\rm H}(E)$ of $\Tan E$ such that
\beq
\Tan E={\rm V}(\pi )\oplus{\rm H}(E)
\label{split}
\eeq
\item
A (global) section of $\pi^1\colon J^1E\to E$;
that is, a mapping $\Psi\colon E\to J^1E$
such that $\pi^1\circ\Psi ={\rm Id}_E$.
\een
\end{prop}
\proof
\quad (1 $\Rightarrow$ 2)\quad
First of all, observe that
$\nabla\colon {\cal X}(E)\to{\cal X}(E)$ is a $\Cinfty (E)$-map
which vanishes when it acts on the vertical vector fields.
Its transposed map is $\nabla^*\colon\df^1(E)\to\df^1(E)$,
which is defined as usually by $\nabla^*\beta :=\beta\circ\nabla$,
for every $\beta\in\df^1(E)$. Moreover, since $\nabla$ is
$\pi$-semibasic, so is $\nabla^*\beta$,
then $\nabla^*(\nabla^*\beta )=\nabla^*\beta$
and hence $\nabla\circ\nabla =\nabla$.
Therefore, $\nabla$ and $\nabla^*$ are projection operators
in ${\cal X}(E)$ and $\df^1(E)$ respectively.
So we have the splittings
$$
{\cal X}(E) = {\rm Im}\nabla \oplus{\rm Ker}\nabla
\quad ; \quad
\df^1(E) = {\rm Im}\nabla^* \oplus{\rm Ker}\nabla^*
$$
with the natural identifications
\beq
({\rm Im}\nabla )'={\rm Ker}\nabla^* \quad ; \quad
({\rm Ker}\nabla )'= {\rm Im}\nabla^*
\label{iden}
\eeq
(where, if $S$ is a submodule of $\vf (E)$, the {\sl incident} or
{\sl annihilator} of $S$ is defined as
$S':=\{\alpha\in\df^1(E)\ | \ \alpha (X)=0 \ , \ \forall X\in S\}$).
Taking this into account, for every $y\in E$,
$\nabla_y\colon\Tan_yE\to\Tan_yE$ induces the splittings
\beq
\Tan_yE = {\rm Im}\nabla_y \oplus{\rm Ker}\nabla_y
\quad ; \quad
\Tan^*_yE = {\rm Im}\nabla^*_y \oplus{\rm Ker}\nabla^*_y
\label{split1}
\eeq
Next we must prove that ${\rm V}_y(\pi )={\rm Ker}\nabla_y$.
But ${\rm V}_y(\pi )\subseteq{\rm Ker}\nabla_y$ and
${\rm Im}\nabla^*_y$ is the set of $\pi$-semibasic forms at $y\in E$,
then we have
${\rm V}_y(\pi )={\rm Ker}\nabla_y$ and, hence,
$$
{\rm Ker}\nabla =\Gamma (E,{\rm V}(\pi ))\equiv\vf^{{\rm V}(\pi)}(E)
\quad ; \quad
{\rm Im}\nabla^* =\Gamma (E,\pi^*\Tan^*M)
$$
So we define
$$
{\rm H}(E):=
\bigcup_{y\in E}\{ \nabla_y(u) \ \vert\ u\in\Tan_yE\}
$$
As a consequence of this, the first of the splittings (\ref{split1})
leads to
\beq
\Tan E={\rm H}(E)\oplus{\rm V}(\pi )
\label{split2}
\eeq
and it allows us to introduce the projections
$$
h\colon\Tan E\longrightarrow{\rm H}(E)\quad ; \quad
v\colon\Tan E\longrightarrow{\rm V}(\pi )
$$
whose transposed maps
$$
h^*\colon{\rm H}^*(\nabla )\longrightarrow\Tan^*E
\quad ; \quad
v^*\colon{\rm V}^*(\pi )\longrightarrow\Tan^*E
$$
are injections which lead to the spliting
\beq
\Tan^*E={\rm H}^*(E)\oplus{\rm V}^*(\pi )
\label{split3}
\eeq
and, taking into account the second equality of
(\ref{split1}) and (\ref{iden}), in a natural way we have the identifications
${\rm H}(E)'$ with ${\rm V}^*(\pi )$ and
${\rm V}(\pi )'$ with ${\rm H}^*(E)$.

\quad\quad (2 $\Rightarrow$ 1)\quad
Given the subbundle ${\rm H}(E)$ and the splitting
$\Tan E={\rm V}(\pi )\oplus{\rm H}(E)$,
let ${\it h}\colon\Tan E\to\Tan E$ and ${\it v}\colon\Tan E\to\Tan E$
be the projections over ${\rm H}(E)$ and ${\rm V}(\pi )$
respectively, which induce the splitting
$X={\it h}(X)+{\it v}(X)$, for every $X\in{\cal X}(E)$.
Then we can define the map
$$
\begin{array}{ccccc}
\nabla&\colon&{\cal X}(E)&\longrightarrow&{\cal X}(E)
\\
& & X & \mapsto & {\it h}(X)
\end{array}
$$
which is a $\Cinfty (E)$-morphism and satisfies
trivially the following properties:
\ben
\item
$\nabla$ vanishes on the vertical vector fields and
therefore $\nabla\in\Gamma (E,\pi^*\Tan^*M)\otimes{\cal X}(E)$.
\item
$\nabla\circ\nabla=\nabla$,
since $\nabla$ is a projection.
\item
if $\alpha\in\Gamma (E,\pi^*\Tan^*M)$
and $X\in{\cal X}(E)$ we have
$$
(\nabla^*\alpha )X=\alpha (\nabla(X))=
\alpha (h(X))=\alpha (h(X)+v(X))=\alpha (X)
$$
because $\alpha$ is semibasic. Therefore $\nabla^*\alpha =\alpha$.
\een

\quad\quad (2 $\Rightarrow$ 3)\quad
Suppose that $\Tan E$ splits into $\Tan E={\rm H}(E)\oplus{\rm V}(\pi)$.
Then there is a natural way of constructing a section of
$\pi^1\colon J^1E\to E$.
In fact, consider $y\in E$ with $\pi(y)=x$, we have
$\Tan_yE={\rm H}_y(E)\oplus{\rm V}_y(\pi)$
and $\Tan_y\pi\vert_{{\rm H}_y(E)}$
is an isomorphism between ${\rm H}_y(E)$ and $\Tan_xM$.
Let $\phi_y\colon U\to E$ be a local section defined in a neigbourhood of
$x$, such that
$$
\phi_y(x)=y \quad  , \quad
\Tan_x\phi_y=(\Tan_y\pi\vert_{{\rm H}_y(E)})^{-1}
$$
then we have a section
$$
\begin{array}{ccccc}
\Psi&\colon&E&\longrightarrow&J^1E
\\
& &y&\mapsto&(j^1\phi_y)(\pi(y))
\end{array}
$$
which is differentiable because the splitting
$\Tan_yE={\rm H}_y(E)\oplus{\rm V}_y(\pi)$
depends differentiabily on $y$.

(3 $\Rightarrow$ 2)\quad
Let $\Psi\colon E\to J^1E$ be a section
and $\bar y\in J^1E$ with $\bar y\map{\pi^1} y\map{\pi} x$.
Observe that $\Psi (y)\in J^1E$ is an equivalence class
of sections $\phi\colon M\to E$ with $\phi (x)=y$,
but the subspace ${\rm Im}\Tan_x\phi$
does not depend on the representative $\phi$, provided it is in this class.
Then, for every $\bar y\in J^1E$ and
$\phi$ a representative of $\bar y=\Psi (y)$, we define
$$
{\rm H}_y(E):={\rm Im}\Tan_x\phi
\quad {\rm and}\quad
{\rm H}(E):=\bigcup_{y\in E}{\rm H}_y(E)
$$
\qed

\begin{definition}
A {\rm connection} in the bundle $\pi\colon E\to M$
is one of the above mentioned equivalent elements.

Then a (global) section $\Psi\colon E\to J^1E$ is said to be a
{\rm 1-jet field} in the bundle $\pi^1\colon J^1E\to E$.
The $\pi$-semibasic form $\nabla$ is called the {\rm connection form}
or {\rm Ehresmann connection}.
The subbundle ${\rm H}(E)$ is called the {\rm horizontal subbundle} of
$\Tan E$ associated with the connection.
The set of sections of ${\rm H}(E)$ will be denoted by ${\cal D}(\Psi )$
and its elements are the {\rm horizontal vector fields}.
\end{definition}

Let $(x^\mu ,y^A)$ be a local system of coordinates in an open set
$U\subset E$. The most general local expression of a
semibasic $1$-form on $E$ with values in $\Tan E$ would be
$$
\nabla=
f_\mu\d x^\mu\otimes\left( g^\nu\derpar{}{x^\nu}+h^A\derpar{}{y^A}\right)
$$
As $\nabla^*$ is the identity on semibasic forms, it follows that
$\nabla^*\d x^\mu =\d x^\mu$,
so the local expression of the connection form $\nabla$ is
$$
\nabla=
\d x^\mu\otimes\left(\derpar{}{x^\mu}+{\mit\Gamma}_\mu^A\derpar{}{y^A}\right)
$$
where ${\mit\Gamma}^A_\mu\in\Cinfty (U)$.
In this system the 1-jet field $\Psi$ is expressed as
$$
\Psi =(x^\mu ,y^A,{\mit\Gamma}_\rho^A(x^\mu ,y^A))
$$

Let $\phi$ be a representative of $\Psi (y)$ with $\phi =(x^\mu ,f^A(x^\mu ))$.
Therefrom $\phi (x)=y$, $\Tan_x\phi=\Psi (y)$ and we have
$$
y=\phi (x)=(x^\mu ,f^A(x^\mu ))=(x^\mu ,y^A)
$$
The matrix of $\Tan_x\phi$ is
\dst
\left(\matrix{{\rm Id}\cr\left(\derpar{f^A}{x^\mu}\right)_x\cr}\right)
\) ,
therefore
\dst\derpar{f^A}{x^\nu}\Big\vert_x={\mit\Gamma}_\nu^A(x^\mu ,y^A)\) .
Now, taking \dst\derpar{}{x^\mu}\Big\vert_x\) as a basis of
$\Tan_xM$, we obtain
$$
{\rm Im}\Tan_x\phi =
\left\{ (\Tan_x\phi )\left(\derpar{}{x^\mu}\right)\right\} =
\left\{\derpar{}{x^\mu}\Big\vert_y+
{\mit\Gamma}_\mu^A(y)\derpar{}{y^A}\Big\vert_y\right\}
$$
hence, ${\rm H}(E)$ is locally spaned by
$$
\left\{\derpar{}{x^\mu}+{\mit\Gamma}_\mu^A(y)\derpar{}{y^A}\right\}
$$

As final remarks, notice that the
splitting (\ref{split2}) induces a further one
$$
{\cal X}(E) = {\rm Im}\nabla \oplus\Gamma (E,{\rm V}(\pi ))
$$
so every vector field $X\in\vf (E)$ splits
into its {\sl horizontal} and {\sl vertical} components:
$$
X=X^{\rm H}+X^{\rm V}=\nabla (X)+(X-\nabla (X))
$$
Locally, this splitting is given by
$$
X = f^\mu\derpar{}{x^\mu}+g^A\derpar{}{y^A} =
f^\mu\left(\derpar{}{x^\mu}+{\mit\Gamma}_\mu^A\derpar{}{y^A}\right)+
\left( g^A-f^\mu{\mit\Gamma}_\mu^A\right)\derpar{}{y^A}
$$
since \dst\derpar{}{x^\mu}+{\mit\Gamma}_\mu^A\derpar{}{y^A}\) and
\dst\derpar{}{y^A}\) generate locally $\Gamma (E,{\rm H}(E))$ and
$\Gamma (E,{\rm V}(\pi ))$, respectivelly.
Observe that, if $X$ is an horizontal vector field,
then $\nabla (X)=X$.

In an analogous way the splitting (\ref{split3}) induces the following one
$$
\df^1(E) = \Gamma (E,\pi^*\Tan^*M)\oplus{\rm Ker}\nabla^* =
\Gamma (E,\pi^*\Tan^*M)\oplus ({\rm Im}\nabla )'
$$
then, for every $\alpha\in\df^1 (E)$, we have
$$
\alpha=\alpha^{\rm H}+\alpha^{\rm B}=
\nabla^*\alpha+(\alpha -\nabla^*\alpha )
$$
whose local expression is
$$
\alpha = F_\mu\d x^\mu+G_A\d y^A =
(F_\mu +G_A{\mit\Gamma}_\mu^A)\d x^\mu+
G_A(\d y^A-{\mit\Gamma}_\mu^A\d x^\mu )
$$
since $\d x^\mu$ and $\d y^A-{\mit\Gamma}_\mu^A\d x^\mu$
generate locally $\Gamma (E,{\rm H}^*(E))$ and
$\Gamma (E,{\rm V}^*(\pi ))$, respectivelly.

\subsection{Other relevant features}

As an interesting remark, we analyze the structure
of the set of connections in $\pi\colon E\to M$.
Then, let $\nabla_1,\nabla_2$ be two connection forms.
The condition $\nabla_1^*\alpha =\nabla_2^*\alpha=0$,
for every semibasic $1$-form $\alpha$, means that
$(\nabla_1-\nabla_2)^*\alpha=0$; that is
$\nabla_1-\nabla_2\in\Gamma (E,\pi^*\Tan^*M)\otimes_E\Gamma (E,{\rm V}(\pi ))$.

However, let $\nabla$ be a connection on $\pi\colon E\to M$ and
$\gamma\in\Gamma (E,\pi^*\Tan^*M)\otimes_E\Gamma (E,{\rm V}(\pi ))$,
then $\nabla +\gamma$ is another connection form. So we have:

\begin{prop}
The set of connection forms on $\pi\colon E\to M$ is an
affine ``space'' over the module of semibasic differential $1$-forms on $E$
with values in ${\rm V}(\pi )$.
\end{prop}

In a local canonical system, if
\dst\nabla =\d x^\mu\otimes
\left(\derpar{}{x^\mu}+{\mit\Gamma}^A_\mu\derpar{}{y^A}\right)\)
and \dst\gamma=\gamma^A_\mu\d x^\mu\otimes\derpar{}{y^A}\)
then
$$
\nabla +\gamma =\d x^\mu\otimes
\left(\derpar{}{x^\mu}+({\mit\Gamma}^A_\mu +\gamma^A_\mu )\derpar{}{y^A}\right)
$$

Finally, another important concept related with connections is:

\begin{definition}
The {\rm curvature} of a connection $\nabla$ is defined as:
$$
{\cal R}(Z_1,Z_2):=
({\rm Id}-\nabla )([\nabla (Z_1),\nabla (Z_2)])=
\inn ([\nabla (Z_1),\nabla (Z_2)])({\rm Id}-\nabla )
$$
for every $Z_1,Z_2\in\vf (M)$.
\label{curva}
\end{definition}

Using the coordinate expressions of the
connection form $\nabla$ and the 1-jet field $\Psi$,
a simple calculation leads to
$$
{\cal R} = \frac{1}{2}
\left(\derpar{{\mit\Gamma}_\eta^B}{x^\mu}-
\derpar{{\mit\Gamma}_\mu^B}{x^\eta}+
{\mit\Gamma}_\mu^A\derpar{{\mit\Gamma}_\eta^B}{y^A}-
{\mit\Gamma}_\eta^A\derpar{{\mit\Gamma}_\mu^B}{y^A}\right)
(\d x^\mu\wedge\d x^\eta )\otimes\derpar{}{y^B}
$$

\subsection{Integrability of jet fields}

\begin{definition}
Let $\Psi\colon E\to J^1E$ be a 1-jet field.
\ben
\item
A section $\phi\colon M\to E$ is said to be an
{\rm integral section} of $\Psi$ iff
$\Psi\circ\phi =j^1\phi $.
\item
$\Psi$ is said to be an {\rm integrable jet field} iff it admits
integral sections.
\een
\end{definition}

One may readily check that, if $(x^\mu ,y^A,v_\mu^A)$ is a natural local
system in $J^1E$ and, in this system,
$\Psi =(x^\mu ,y^A,{\mit\Gamma}_\rho^A(x^\mu ,y^A))$ and
$\phi =(x^\mu ,f^A(x^\nu ))$,
then $\phi$ is an integral section of $\Psi$ if, and only if,
$\phi$ is a solution of the following system of
partial differential equations
\beq
\derpar{f^A}{x^\mu}={\mit\Gamma}_\mu^A\circ\phi
\label{condin}
\eeq

The integrable jet fields can be characterized as follows:

\begin{prop}
The following assertions on a 1-jet field $\Psi$ are equivalent:
\ben
\item
The 1-jet field $\Psi$ is integrable.
\item
The curvature of the connection form $\nabla$
associated with $\Psi$ is nule.
\item
${\cal D}(\Psi )$ is an involutive distribution.
\een
\end{prop}
( {\sl Proof} )\quad
(1 $\Leftrightarrow$ 2)\quad
Notice that if $\phi$ is an integral section of $\Psi$, then the distribution
${\cal D}(\Psi )$ is tangent to the image of $\Psi$ and conversely.

\quad\quad (2 $\Leftrightarrow$ 3)\quad
 From the definition (\ref{curva}) we obtain that, if ${\cal R}=0$, then
$$
\nabla([\nabla (Z_1),\nabla (Z_2)])=[\nabla (Z_1),\nabla (Z_2)]
$$
hence, the horizontal distribution ${\cal D}(\Psi )$ is involutive.

Conversely, if ${\cal D}(\Psi )$ is involutive,
as $\nabla$ is the identity on ${\cal D}(\Psi )$,
the last equation follows, so ${\cal R}=0$.
\qed

{\bf Remark}:
\bit
\item
According to this proposition,
from the local expression of ${\cal R}$ we obtain the local
integrability conditions of the equations (\ref{condin}).
\eit

\section{Glossary of notation}

\begin{tabular}{ll}
$\Cinfty (A)$ &
Smooth functions in the manifold $A$.
\\
$\vf (A)$ &
Vector fields in the manifold $A$.
\\
$\df^p(A)$ &
Differential $p$-forms in the manifold $A$.
\\
$A'$ &
Incident set to $A$.
\\
${\rm Hom}(A,B)$ &
Homomorphisms from $A$ to $B$.
\\
$\inn (X)\alpha$ &
Contraction of the elements $X$ and $\alpha$.
\\
$\Lie (X)\alpha$ &
Lie derivative of the element $\alpha$ along the vector field $X$.
\\
$D_v\alpha$ &
Directional derivative of $\alpha$ along the direction of $v$.
\\
$\pi\colon E\to M$ &
Differentiable fiber bundle over a (orientable) manifold $M$.
\\
$\Gamma (A,B)$ , $\Gamma ({\rm pr})$ &
Sections of the bundle $pr\colon B\to A$.
\\
$\Gamma_U(A,B)$ &
Local sections (in an open set $U\subset A$).
\\
$\Gamma_c(A,B)$ &
Compact supported sections .
\\
$\df^p(A,B)=\df^p(A)\otimes_A\Gamma (B)$ &
Differential $p$-forms in the manifold $A$ with values in the bundle $B$.
\\
$\pi^1\colon J^1E\to E$ &
$1$-jet bundle.
\\
$\bar\pi^1\colon J^1E\to M$ &
Fiber bundle induced by $\bar\pi^1=\pi\circ\pi^1$.
\\
$J^1_yE$ &
Fiber of $J^1E$ at $y\in E$.
\\
$\bar y\in J^1E$ &
$1$-jet (element of $J^1E$).
\\
$(x^\mu ,y^A,v^A_\mu )$ &
Local coordinate system in $J^1E$.
\\
${\rm V}(pr)$ &
Vertical subbundle with respect to the projection $pr\colon A\to M$.
\\
${\cal X}^{{\rm V}(pr)}(B)$ &
Vertical vector fields in $B$ with respect to the projection
$pr\colon B\to A$.
\\
$\phi\colon M\to E$ &
Section of $\pi$.
\\
$\psi\colon M\to J^1E$ &
Section of $\bar\pi^1$.
\\
$\d^v_y\phi$ &
Vertical differential of $\phi$ at $y\in E$.
\\
$\theta$ &
Structure canonical form of $J^1E$.
\\
${\cal M}_c$ &
Contact module or Cartan distribution.
\\
$\theta^A$ &
Element of a local basis of ${\cal M}_c$.
\\
${\cal I}({\cal M}_c)$ &
Ideal generated by the contact module.
\\
${\cal V}$ , ${\cal S}$ &
Vertical endomorphisms in $J^1E$.
\\
$j^1\phi$, ($j^\psi$) &
Canonical prolongation of $\phi$ ($\psi$) to $J^1E$ ($J^1J^1E$).
\\
$j^1\Phi$ &
Canonical prolongation of a diffeomorphism $\Phi\colon E\to E$ to $J^1E$.
\\
$j^1Z$ &
Canonical prolongation of a vector field $Z\in\vf (E)$ to $J^1E$.
\\
$\nabla$ &
Connection form (Ehresmann connection).
\\
${\rm H}(E)$ &
Horizontal subbundle (induced by a connection).
\\
$\Psi\colon E\to J^E$ &
(Global) section of $\pi^1$. ($1$-jet field in $J^1E$).
\\
${\cal D}({\Psi})$ &
Set of sections of ${\rm H}(E)$ (horizontal vector fields).
\\
${\cal R}$ &
Curvature of a connection.
\\
$\Lag$ &
Lagrangian density.
\\
$\lag$ &
Lagrangian function.
\\
${\bf L}$ &
Lagrangian functional.
\\
$\vartheta_{\Lag}$ &
Lagrangian canonical form.
\\
$\Theta_{\Lag}$ &
Poincar\'e-Cartan $(n+1)$-form.
\\
$\Omega_{\Lag}$ &
Poincar\'e-Cartan $(n+2)$-form.
\\
$\ls$ &
Lagrangian system
\\
${\cal E}^\nabla_{\Lag}$ &
Density of lagrangian energy associated with $\Lag$ and $\nabla$.
\\
${\rm E}^{\nabla}_{\Lag}$ &
Energy associated with $\Lag$ and $\nabla$.
\\
$\pi^1_1\colon J^1J^1E\to J^1E$ &
Repeated $1$-jet bundle.
\\
${\bf y}\in J^1J^1E$ &
Element of $J^1J^1E$.
\\
${\cal Y}\colon J^1E\to J^1J^1E$ &
(Global) section of $\pi^1_1$ ($1$-jet field in $J^1J^1E$).
\\
$\delta_Z\lag$ &
Total variation of $\lag$ (under the action of $Z\in\vf (E)$).
\\
$(j^1\phi )^*{\rm J}Z:=(j^1\phi )^*(\inn (j^1Z)\Theta_{\Lag})$ &
Noether's current associated with $Z$ (for a critical section $\phi$).
\\
\end{tabular}

\subsection*{Acknowledgments}

We thank Dr. X. Gr\`acia-Sabat\'e (U.P.C.) for discussions
concerning some subjects on Differential Geometry and
Dr. L.A. Ibort (U.C.M.) for his comments and suggestions.
We also thank Mr. Jeff Palmer for his assistance in
preparing the English version of the manuscript.

\end{document}